\documentclass[12pt,english]{amsart}
\usepackage{palatino}
\usepackage[T1]{fontenc}
\usepackage[latin1]{inputenc}
\usepackage{geometry}
\usepackage{pifont}
\usepackage{color}
\usepackage{amssymb}

\makeatletter

 \theoremstyle{plain}
 \theoremstyle{remark}    
 \newtheorem*{acknowledgement*}{Acknowledgement} 
 \theoremstyle{definition}
 \newtheorem{defn}{Definition}
 \theoremstyle{plain}    
 \newtheorem{lem}{Lemma} 
 \theoremstyle{remark}
 \newtheorem{rem}{Remark}
 \theoremstyle{plain}    
 \newtheorem{thm}{Theorem} 
 \theoremstyle{plain}    
 \newtheorem{prop}{Proposition} 
 \theoremstyle{remark}    
 \newtheorem{claim}{Claim}
 \theoremstyle{plain}    
 \newtheorem{cor}{Corollary} 
 \theoremstyle{definition}
  \newtheorem{example}{Example}

\usepackage{euscript}
\usepackage{amscd}

\newcommand{\Out}{\operatorname{Out}}
\newcommand{\Aut}{\operatorname{Aut}}
\newcommand{\Pic}{\operatorname{Pic}}
\newcommand{\StatPic}{\operatorname{StatPic}}
\newcommand{\Hol}{\operatorname{Hol}}
\newcommand{\Mod}{\operatorname{Mod}}
\newcommand{\PMod}{\operatorname{PMod}}
\newcommand{\Diff}{\operatorname{Diff}}
\newcommand{\id}{\operatorname{id}}
\newcommand{\pr}{\operatorname{pr}}
\newcommand{\Symp}{\operatorname{Symp}}
\newcommand{\Poiss}{\operatorname{Poiss}}
\newcommand{\op}{\operatorname{op}}
\newcommand{\Bis}{\operatorname{Bis}}

\newcommand{\fix}{\operatorname{fix}}

\newcommand{\End}{\operatorname{End}}

\let\mathcal\EuScript
\AtBeginDocument{
  
}

\usepackage{babel}
\makeatother
\begin{document}

\title{Picard groups of topologically stable Poisson structures }

\author{Olga Radko{*} and  Dimitri Shlyakhtenko{*}{*}}

\thanks{{*}Research supported by a VIGRE Postdoctoral Fellowship. \\
{*}{*}Research supported by a Sloan Foundation Fellowship and NSF
grant DMS-0102332.}

\address{Department of Mathematics, UCLA, Los Angeles, CA, 90095, USA}

\begin{abstract}
We compute the group of Morita self-equivalences (the Picard group)
of a Poisson structure on an orientable surface, under the assumption
that the degeneracies of the Poisson tensor are linear. The answer
involves mapping class groups of surfaces, i.e., groups of isotopy
classes of diffeomorphisms. We also show that the Picard group of
these structures coincides with the group of outer Poisson automorphisms.
\end{abstract}
\maketitle

\section{Introduction\label{sec:Introduction}}

There are many similarities between Poisson geometry and the theory
of associative algebras (see e.g. the book \cite{CannasdaSilva-Weinstein:GometricModelsNCAlg}).
Based on the notion of Morita equivalence for Poisson manifolds introduced
by Xu \cite{Xu-ME-Poisson}, Bursztyn and Weinstein \cite{BW1} have
recently defined the \emph{Picard group} $\Pic(P)$ of an integrable
Poisson manifold $P$ as the group of all Morita equivalences between
$P$ and itself (see Section \ref{sub:Picard_group} for a definition).
The group operation is the relative tensor product $\otimes_{P}$
of bimodules (see \cite{Xu-ME-symp_groupoids} ), and the identity
bimodule is the source simply connected symplectic groupoid $\Gamma(P)$.
The Picard group contains the group of outer Poisson automorphisms
of $P$ in a natural way, but it can in principle be strictly larger
(such is the case, e.g., for an open symplectic surface; Section \ref{sec:trivializations}). 

The aim of this paper is to present a complete computation of the
Picard group for a certain class of Poisson structures on compact
connected oriented surfaces. We consider Poisson structures (called
\emph{topologically stable structures}, or TSS, for short) which are
non-degenerate almost everywhere on the surface, except that they
have linear degeneracies on a finite set of simple closed curves.
Although TSS are sufficiently generic (the set of these structures
is open and dense in the vector space of all Poisson structures on
a given surface), in many ways they resemble the symplectic structures.
In \cite{Radko-classification}, the first author has obtained a complete
description of the moduli space of isomorphism classes of TSS by giving
a complete set of explicit invariants. A complete set of criteria
for Morita equivalence of two TSS was found in \cite{BW1}.

The problem of computing the Picard group of a TSS was posed in \cite{BW1}.
Our main result is the following. Let $\pi$ be a TSS on a surface
$P$. Let $Z$ be the zero set of this structure, consisting of a
disjoint union of a finite number of simple closed curves $T_{1},T_{2},\ldots T_{n}$
on $P$. A modular vector field (\cite{Weinstein_modular}) of $(P,\pi)$
restricts to a canonical nowhere zero vector field $\xi$ along the
zero set $Z$. Consider the group $D(\pi)$ of all diffeomorphisms
$\phi:P\to P$, so that

\begin{itemize}
\item $\phi$ preserves the zero set, $\phi(Z)=Z$;
\item $\phi$ preserves the restriction of a modular vector field, $\phi_{*}\xi=\xi$.
\end{itemize}
Let $D_{0}(\pi)$ be the subgroup of $D(\pi)$ consisting of diffeomorphisms
fixing a neighborhood of $Z$ pointwise. (This actually implies that
$D_{0}(\pi)$ is precisely the subgroup of diffeomorphisms which preserve
the leaf space of $(P,\pi)$ pointwise). Then \begin{equation}
\Pic(P,\pi)\cong D(\pi)/\textrm{isotopy by elements of $D_{0}(\pi)$}.\label{eq:Pic_as_equiv_classes}\end{equation}
In other words, the Picard group can be identified with the set of
equivalence classes of diffeomorphisms in $D(\pi)$, two diffeomorphisms
$\phi_{1}$ and $\phi_{2}$ being equivalent if $\phi_{1}\circ\phi_{2}^{-1}$
can be connected to the identity by a continuous path of diffeomorphisms
lying in $D_{0}(\pi)$.

Our description of the Picard group of a TSS gives also an explicit
formula for it, involving the mapping class groups. If $Z=\bigsqcup_{i\in I}T_{i}$
is the zero set of $\pi$, and $P\setminus Z=\bigsqcup_{j\in J}L_{j}$
is the decomposition of the complement of the zero set into the $2$-dimensional
symplectic leaves, then\begin{equation}
\Pic(P,\pi)\cong(\prod_{i\in I}\mathbb{T}\times\prod_{j\in J}\mathcal{M}(L_{j}))\rtimes G.\label{eq:Pic_as_semi_dir_product}\end{equation}
Here for each zero curve, $\mathbb{T}$ is the $1$-torus of translations
by the flow of the restriction of a modular vector field to this curve.
$\mathcal{M}(S)$ stands for the mapping class group of an open surface
$S$:\[
\mathcal{M}(S)=\{\textrm{diffeomorphisms of }S\textrm{ fixing pointwise a neighborhood of infinity}\}/\textrm{isotopy}.\]
 Mapping class groups of surfaces are well understood, and can be
explicitly described in terms of generators and relations. However,
their appearance in relation to Picard groups of TSS was quite unexpected
and surprising. 

The group $G$ in (\ref{eq:Pic_as_semi_dir_product}) is the (discrete)
group of automorphisms of a labeled graph encoding the Morita equivalence
class of $(P,\pi)$. Automorphisms in $G$ permute the leaves, and
thus act naturally on the set $I$ of zero curves and the set $J$
of $2$-dimensional leaves. 

We finish the introduction with an outline of the proof of the main
result. In order to prove the isomorphism (\ref{eq:Pic_as_equiv_classes}),
one needs to be able to prove that:

\begin{enumerate}
\item Every class of an element of $D(\pi)$ gives rise to a Morita self-equivalence
bimodule;
\item Every Morita self-equivalence bimodule arises in this way. 
\end{enumerate}
The first part is relatively easy: if $\phi\in D(\pi)$, then a Moser-type
argument (like in \cite{Radko-classification}) implies that $\phi$
is isotopic to a Poisson automorphism $\phi_{0}$. Then one simply
assigns to the isotopy class of $\phi$ the bimodule $\Gamma_{\phi_{0}}(P)$
associated to this Poisson automorphism (see Section \ref{sub:Bimodules_from_autos}
for a definition). This results in a map $j:(D(\pi)/\textrm{isotopy})\to\Pic(P,\pi)$.

The second part is more difficult. We need to show that the map $j$
is onto. Clearly, it is sufficient to show that if $X$ is a Morita
self-equivalence bimodule over $P$, then for some bimodule $X'$
in the image of $j$, the relative tensor product $X\otimes_{P}X'$
is isomorphic to the identity bimodule. In particular, it is necessary
to be able to determine when a Morita equivalence bimodule is trivial
(i.e., isomorphic to the identity bimodule). 

In algebra, one meets a very simple criterion of triviality: an $A,A$-bimodule
$X$ is trivial if and only if it contains a vector $v_{0}$ which
is bicyclic (i.e., $X=Av_{0}A$) and central (i.e., $av_{0}=v_{0}a$
$\forall a\in A$). In that case, the map $a\mapsto av_{0}$ gives
an isomorphism between the identity bimodule $A$ and $X$. In Poisson
geometry the analogous notion is that of an \emph{identity bisection}
(Definition \ref{def:identity_bisection}). The existence of an identity
bisection gives a criterion for a bimodule over an integrable Poisson
manifold to be trivial.

There are three obstructions to existence of an identity bisection
for a Morita self-equivalence bimodule $X$ over $P$:

\begin{description}
\item [O1]$X$ may induce a non-trivial automorphism of the leaf space
(see Section \ref{sub:The-Static-Picard}). 
\item [O2]The restriction of $X$ to the open symplectic manifold $P\setminus Z$
may be non-trivial in the Picard group of $P\setminus Z$.
\item [O3]Even if the restriction of $X$ to $P\setminus Z$ is trivial
and thus has an identity bisection defined on $P\setminus Z$, this
bisection may fail to extend to all of $P$.
\end{description}
The first obstruction vanishes iff $X$ is a \emph{static} $P,P$-bimodule
(i.e., it induces the identity map on the leaf space of $P$). One
can easily check that for any bimodule $X\in\Pic(P,\pi)$, there is
a bimodule $X'$ in the image of $j$, so that $X\otimes_{P}X'$ is
static. Thus we can restrict our consideration to static bimodules.

If $X$ is static, then we show that the restriction of $X$ to a
cylindrical neighborhood of a zero curve automatically has exactly
one identity bisection on this neighborhood. This gives a \emph{local
identity bisection} (see Definition \ref{def:local_identity_bisection})
of $X$ in a neighborhood of the zero set $Z$. The existence and
uniqueness of this local identity bisection has a lot to do with the
presence of a very special local symmetry of the Poisson structure,
given by the flow of the modular vector field \cite{Weinstein_modular},
and the fact that modular vector fields are {}``preserved'' by Morita
equivalences \cite{Ginzburg-Lu}.

Fix a $2$-dimensional symplectic leaf $L\subset P\setminus Z$. By
a result of Bursztyn and Weinstein \cite{BW1}, $\Pic(L)\cong\Out(\pi_{1}(L))$.
 However, since any $X\in\Pic(P)$ has a local identity bisection
near the zero set, the restriction of $X$ to $L$ must be {}``trivial
near the boundary'' of $L$. Thus it cannot give an arbitrary element
in $\Out(\pi_{1}(L))$; it must give elements that {}``preserve peripheral
structure'' of $L$. (See Section \ref{sec:trivializations}). 

By the Dehn-Nielsen-Baer's theorem (see Theorem \ref{thm:Dehn-Nielsen-Baer}),
one can find a symplectomorphism defined on $L$ and trivial near
its boundary, which induces any element of $\Out(\pi_{1}(L))$, preserving
peripheral structure. Thus given a bimodule $X\in\Pic(P,\pi)$, with
vanishing O1, one can use the Dehn-Nielsen-Baer's theorem to find
a bimodule $X'$ in the image of $j$, so that $X\otimes_{P}X'$ has
the property that both obstructions O1 and O2 vanish. 

Assuming that obstructions O1 and O2 vanish, the last obstruction
O3 is related to the behavior of the identity bisection(s) of the
restriction of $X$ to $P\setminus Z$ near $Z$. The issue is whether
the identity bisection defined on a leaf $L\subset P\setminus Z$
matches the unique local identity bisection defined in a neighborhood
of $\partial L\subset Z$. It turns out that the obstruction is integer-valued
in a neighborhood of each zero curve. Moreover, the value of the obstruction
associated to the Dehn twist near a zero curve is exactly $1$. This
implies that given $X\in\Pic(P,\pi)$, with O1 and O2 vanishing, it
is possible to find $X'$ associated to a combination of Dehn twists
(and thus in the image of $j$), so that all obstructions O1, O2 and
O3 vanish for $X\otimes_{P}X'$. These Dehn twists are special elements
of mapping class groups $\mathcal{M}(L_{j})$ (see (\ref{eq:Pic_as_semi_dir_product})).

\begin{acknowledgement*}
We are very grateful to Henrique Bursztyn and Alan Weinstein for the
discussions of this work and for their useful comments. We would also
like to acknowledge our conversations related to mapping class groups
with Alex Bene and Jeff Mess.
\end{acknowledgement*}

\section{Preliminaries on Morita equivalence and Picard groups\label{sec:preliminaries}}

\subsection{Symplectic Groupoids and Modules over Poisson manifolds.}

\subsubsection{Symplectic Groupoids.}

Let $(P,\pi)$ be an integrable Poisson manifold and let $\Gamma(P)$
be its (source-connected and simply connected) \emph{symplectic groupoid}.
The elements of $\Gamma(P)$ can be thought of as classes of \emph{cotangent
paths} up to cotangent homotopy%
\footnote{The notion of cotangent homotopy is somewhat technical and is not
explicitly needed in the present paper. We refer the reader to \cite{CF2}
for details.%
} (see \cite{CF1,CF2}). 

A cotangent path is a pair $(a,\gamma)$, where $a:\,[0,1]\to T^{*}P$
is a path in the cotangent bundle and $\gamma$ is a path $\gamma:\,[0,1]\to P$
on the manifold (called the base path of $a$) such that

\begin{itemize}
\item $\textrm{pr}(a(t))=\gamma(t)$;
\item $\frac{d\gamma}{dt}=\tilde{\pi}(a(t))$.
\end{itemize}
Here $\textrm{pr}:T^{*}P\to P$ is the natural projection and $\tilde{\pi}:\  T^{*}P\to TP$
is the bundle map associated to the Poisson structure. The symplectic
structure on $\Gamma(P)$ is obtained from the natural symplectic
structure on its Lie algebroid, $T^{*}P$. The source and target maps
$s,t:\ \Gamma(P)\to P$ send a cotangent path $a$ over a base path
$\gamma$ to the beginning and end of $\gamma$, respectively.

\subsubsection{The symplectic groupoid of a symplectic manifold}

In the special case of a symplectic manifold $S$, the bundle map
$\tilde{\pi}:\  T^{*}S\to TS$ is invertible and therefore a cotangent
path is uniquely determined by its base path. In this case, cotangent
homotopy is equivalent to the usual homotopy of base paths. In particular,
as a groupoid, $\Gamma(S)$ can be identified with the fundamental
groupoid of $S$; in other words, an element in $\Gamma(S)$ can be
viewed as a class of a (not necessarily closed) path in $S$, considered
up to a homotopy fixing the endpoints. The natural symplectic structure
on $\Gamma(S)$ is such that $\Gamma(S)=\frac{\tilde{S}\times\tilde{S}^{\op}}{\pi_{1}(S)}$,
where the universal cover $\tilde{S}$ is taken with the pull-back
of the original symplectic structure on $S$, and $\tilde{S}^{\op}$
denotes the same manifold with the negative of the symplectic structure.
The \emph{isotropy group} $\Gamma_{p}(S)=s^{-1}(p)\cap t^{-1}(p)$
of $\Gamma(S)$ at a point $p\in S$ is canonically isomorphic to
the fundamental group $\pi_{1}(S,p)$ of $S$ at $p$, and is therefore
discrete.

\subsubsection{Modules over Poisson manifolds.}

If one thinks of Poisson manifolds as semi-classical analogs of associative
algebras, a symplectic manifold $(X,\Omega)$ should be thought of
as an analog of the algebra $\End(V)$, where $V$ is a vector space.
The structure of a left (respectively, right) module over an algebra
$A$ on a vector space $V$ is given by an algebra (anti) homomorphism
from $A$ to $\End(V)$. 

The analogous notion in Poisson geometry is that of a \emph{left (right)
module} over a Poisson manifold $P$. This is defined as a complete
(anti)-symplectic realization $J:\ (X,\Omega)\to(P,\pi)$, i.e., an
(anti)-Poisson map from a symplectic manifold $(X,\Omega)$ to $(P,\pi)$. 

Any symplectic realization $J:\ (X,\Omega)\to(P,\pi)$ induces a canonical
action of the Lie algebroid $T^{*}P$ of the Poisson manifold on $X$.
In other words, there is a Lie algebra homomorphism $\Phi_{X}:\Gamma(T^{*}P)\to\Gamma(TX)$
given by $\Phi_{X}(\alpha)=\tilde{\Omega}^{-1}(J^{*}\alpha)$ for
$\alpha\in\Omega^{1}(P)=\Gamma(T^{*}P)$. If $P$ is integrable and
the symplectic realization $J:X\to P$ is complete, this Lie algebroid
action integrates to an action of the symplectic groupoid $\Gamma(P)$
on $X$. In this case, the map $J$ is called the moment map for the
groupoid action. Thus, there is a correspondence between the modules
over an integrable Poisson manifold and the actions of its symplectic
groupoid.

\subsection{Bimodules and Morita equivalence.}

Most of the definitions and constructions in this section can be found
in \cite{BW1}. We refer the reader to that paper for more details
on Picard groups of Poisson manifolds and symplectic groupoids.

\subsubsection{Bimodules over Poisson manifolds.}

For two Poisson manifolds, a $P_{1},P_{2}$-\emph{bimodule} is a symplectic
manifold $X$ which is a left $P_{1}$-module and a right $P_{2}$-module,
such that the corresponding (left and right) Lie algebroid (or, equivalently,
groupoid) actions commute. Thus, we have a pair $P_{1}\stackrel{s_{X}}{\leftarrow}X\stackrel{t_{X}}{\rightarrow}P_{2}$
of maps such that

\begin{enumerate}
\item $s_{X}$ is a complete Poisson map and $t_{X}$ is a complete anti-Poisson
map;
\item $s_{X}$ and $t_{X}$ are surjective submersions with connected and
simply-connected fibers;
\item $\{ s_{X}^{*}C^{\infty}(P_{1}),t_{X}^{*}C^{\infty}(P_{2})\}=0$; 
\end{enumerate}
The last condition is equivalent to the requirement the $s_{X}$-
and $t_{X}$- fibers are symplectically orthogonal to each other in
$X$. We will frequently refer to $X$ as a $P_{1},P_{2}$-bimodule,
implicitly denoting by $s_{X}$ and $t_{X}$ the associated bimodule
maps.

Two bimodules are isomorphic if there is a symplectomorphism preserving
the bimodule structure.

\subsubsection{The relative tensor product $\otimes_{P}$}

Given a $P_{1},P_{2}$-bimodule $X$ and a $P_{2},P_{3}$-bimodule
$Y$, one can define their relative tensor product $X\otimes_{P_{2}}Y$
(see \cite{HilsumSkandalis} for the case of bimodules over groupoids,
and \cite{Xu-ME-symp_groupoids} for the Poisson case) as the orbit
space\[
X\otimes_{P_{2}}Y=X\times_{P_{2}}Y/\Gamma(P_{2}).\]
Here $X\times_{P_{2}}Y=\{(x,y)\in X\times Y\ :\  t_{X}(x)=s_{Y}(y)\}$
is the fibered product, and the action of the symplectic groupoid
$\Gamma(P_{2})$ is given by $g\cdot(x,y)=(xg,\  g^{-1}y)\ $ for
$g\in\Gamma(P_{2})$, $(x,y)\in X\times_{P_{2}}Y$. 

If $X\otimes_{P_{2}}Y$ is a smooth manifold, it is automatically
a $P_{1},P_{3}$-bimodule. The symplectic structure of $X\otimes_{P_{2}}Y$
was first described in \cite{Xu-ME-symp_groupoids}. To ensure the
smoothness, it is enough to assume that $Y$ is a \emph{left principal
bimodule} (cf. \cite{Landsman_tensor_products,BW1}), i.e., that the
action of $\Gamma(P_{1})$ is free and transitive on the fibers of
$s_{Y}$. 

For a Poisson manifold $P$, the symplectic groupoid $\Gamma(P)$,
considered as a $P,P$-bimodule using its source and target as the
bimodule maps, plays the role of the identity bimodule: for a right
$P$-module $X$ and a left $P$-module $Y$, we have\[
X\otimes_{P}\Gamma(P)\cong X,\qquad\Gamma(P)\otimes_{P}Y\cong Y.\]

\subsection{The Picard group.\label{sub:Picard_group}}

Two (integrable) Poisson manifolds $P_{1}$ and $P_{2}$ are called
\emph{Morita equivalent} if there exists a $P_{1},P_{2}$-bimodule
$X$ with the property that it is \emph{invertible}, i.e., there is
a $P_{2},P_{1}$-bimodule $Y$ so that $X\otimes_{P_{2}}Y\cong\Gamma(P_{1})$
and $Y\otimes_{P_{1}}X\cong\Gamma(P_{2})$. By a result in \cite{BW1},
invertibility is equivalent to the requirement that the actions of
$\Gamma(P_{1})$ and $\Gamma(P_{2})$ on the fibers of $t_{X}$ and
$s_{X}$ are free and transitive. If $X$ is invertible, the inverse
is necessarily isomorphic to the \emph{opposite bimodule} $X^{\op}$
(obtained by switching the bimodule maps and changing the symplectic
structure of $X$ to its negative).

The \emph{Picard group} of a Poisson manifold $P$ is the group $\Pic(P)$
of all isomorphism classes of Morita equivalence $P,P$-bimodules,
considered with the operation of the relative tensor product. In analogy
with the algebraic case, the Picard group can be considered as a group
of generalized automorphisms of the structure.

\subsubsection{Bimodules associated to Poisson isomorphisms.\label{sub:Bimodules_from_autos}}

Each Poisson diffeomorphism $\varphi\in\Poiss(P)$ gives rise to a
self-equivalence bimodule $\Gamma_{\varphi}(P)\in\textrm{Pic}(P)$.
As a symplectic manifold, $\Gamma_{\phi}(P)=\Gamma(P)$. The bimodule
structure of $\Gamma_{\varphi}(P)$ is given by\[
P\stackrel{s}{\leftarrow}\Gamma(P)\stackrel{\varphi\circ t}{\longrightarrow}P.\]
In other words, the bimodule structure of $\Gamma_{\phi}(P)$ is obtained
from that of $\Gamma(P)$ by twisting the target map $t$ by $\phi$. 

It follows that there is a group homomorphism\[
j:\Poiss(P)\to\Pic(P)\]
given by $j(\varphi)=\Gamma_{\varphi}(P).$ It turns out \cite{BW1}
that the kernel of this homomorphism consists exactly of inner Poisson
automorphisms, i.e., automorphisms that can be induced by the action
of the group of bisections of $\Gamma(P)$.

\subsubsection{The Static Picard group.\label{sub:The-Static-Picard}}

Let $\mathcal{L}(P)$ be the leaf space of $P$, regarded as a topological
space with the quotient topology. Let $\textrm{Aut}(\mathcal{L}(P))$
be the group of its homeomorphisms. 

A $P,P$ -bimodule $P\stackrel{s_{X}}{\leftarrow}X\stackrel{t_{X}}{\rightarrow}P$
defines a homeomorphism $\phi_{X}\in\textrm{Aut}(\mathcal{L}(P))$
by $\phi_{X}:L\mapsto t_{X}\circ s_{X}^{-1}(L)$. This gives a group
homomorphism $h:\Pic(P)\to\Aut(\mathcal{L}(P))$. Its kernel forms
a subgroup $\StatPic(P)\subseteq\Pic(P)$, called the \emph{static
Picard group} of $P$, which consists of self-equivalence bimodules
fixing the leaf space pointwise. The computations of the Picard group
can sometimes be simplified by using the following exact sequence:\begin{equation}
1\to\StatPic(P)\to\Pic(P)\stackrel{h}{\to}\Aut(\mathcal{L}(P)).\label{eq:Picard_exact_sequence}\end{equation}

\subsection{Identity Bisections.\label{sub:Identity-Bisections.}}

Let $X$ be a bimodule over a Poisson manifold $P$. 

\begin{defn}
\label{def:identity_bisection}A map $\varepsilon:P\to X$ is called
an \emph{identity bisection} of $X$ if the following conditions are
satisfied:
\end{defn}
\begin{enumerate}
\item $s_{X}\circ\varepsilon=t_{X}\circ\varepsilon=\id\in\Diff(P)$;
\item for any $\gamma\in\Gamma(P)$ with $s(\gamma)=p$, $t(\gamma)=q$,
one has\[
\gamma\cdot\varepsilon(q)=\varepsilon(p)\cdot\gamma,\]
i.e. the actions of the symplectic groupoid commute with taking the
bisection. 
\item $\varepsilon(P)\subset X$ is a Lagrangian submanifold.
\end{enumerate}
Conditions 1 and 2 are actually equivalent for a bisection $\varepsilon:P\to X$
defined on all of $P$ (see \cite{Coste-Dazord-A.W.}). However, we
prefer to distinguish between them, since locally these conditions
are not equivalent. Indeed, if $U\subset P$ is any subset that intersects
every leaf of $P$, and $\varepsilon:U\to X|_{U}$ satisfies Condition
2 (i.e., for all $\gamma\in\Gamma(P)$ with $s(\gamma)=p\in U$, $t(\gamma)=q\in U$,
one has $\gamma\cdot\varepsilon(q)=\varepsilon(p)\cdot\gamma$), then
$\varepsilon$ can be uniquely extended to all of $P$. On the other
hand, Condition 1 may be satisfied on $U$ without there being an
identity bisection defined on all of $P$. For this reason we make
the following definition:

\begin{defn}
\label{def:local_identity_bisection}Let $U\subset P$ be a subset
and let $X$ be a $P,P$-bimodule. We say that $\varepsilon:U\to X|_{U}$
is a \emph{local identity bisection}, if
\begin{enumerate}
\item $s_{X}\circ\varepsilon(p)=t_{X}\circ\varepsilon(p)=p$ for all $p\in U$;
\item $\varepsilon(U)\subset X$ is a Lagrangian submanifold.
\end{enumerate}
\end{defn}
\begin{lem}
\label{lem:identityBisectionTrivial}A Morita self-equivalence $P,P$-bimodule
$X$ is isomorphic to the identity bimodule $\Gamma(P)$ if and only
if $X$ has an identity bisection. 
\end{lem}
\begin{proof}
Assume that $X$ has an identity bisection $\varepsilon:P\to X$ .
Then\[
\gamma\mapsto\gamma\cdot\varepsilon(t(\gamma))=\varepsilon(s(\gamma))\cdot\gamma\]
is a bimodule map from $\Gamma(P)$ to $X$. This map is injective
and surjective, since the action of $\Gamma(P)$ on $X$ is free and
transitive on the fibers. Thus, we obtain a bimodule isomorphism.
The condition that $\varepsilon(P)$ is Lagrangian guarantees that
this map is a symplectomorphism. Conversely, if $X$ is isomorphic
to $\Gamma(P)$, the image of the identity bisection of $\Gamma(P)$
under this isomorphism is an identity bisection in $X$. 
\end{proof}

\section{Mapping Class Groups of surfaces\label{sec:Mapping-Class-Groups} }

\subsection{The groups $\Mod(S)$, $\PMod(S)$ and $\mathcal{M}(S)$.}

Throughout this section, unless explicitly stated otherwise, let $S$
be a connected oriented surface. For simplicity, we will assume that
$S$ is the interior of a surface $\bar{S}$ with boundary $\partial S$.
We will be mostly concerned with the case that $\partial S$ is non-empty.

We denote by $\Diff(S)$ the group of diffeomorphisms of $S$, and
by $\Diff(S\fix\partial S)$ its subgroup of diffeomorphisms preserving
pointwise a neighborhood of the boundary of $S$ (note that the diffeomorphisms
in $\Diff(S\fix\partial S)$ are automatically orientation-preserving
as long as $\partial S\neq\emptyset$). Both $\Diff(S)$ and $\Diff(S\fix\partial S)$
are topological groups with the $C^{\infty}$-topology. Mapping class
groups arise from the consideration of the connected components of
these groups. Following the notation of the survey \cite{MappingClassGroups},
we consider:

\begin{enumerate}
\item \emph{The mapping class group} $\Mod(S)=\pi_{0}(\Diff(S))$, defined
as the group of isotopy classes of all diffeomorphisms of $S$.
\item \emph{The pure mapping class group} $\PMod(S)$, defined as the subgroup
of $\Mod(S)$ generated by diffeomorphisms that lie in $\Diff(S\fix\partial S)$.
\item $\mathcal{M}(S)=\pi_{0}(\Diff(S\fix\partial S))$. 
\end{enumerate}
\begin{rem}
Note that in the case of closed surfaces (i.e., $\partial S=\emptyset$),
we have $\mathcal{M}(S)=\PMod(S)=\Mod(S)$. 
\end{rem}
Let $\varphi\in\Diff(S)$ be a diffeomorphism and $\varphi_{*}:\pi_{1}(S,p)\to\pi_{1}(S,\varphi(p))$
be the induced group homomorphism. This gives an automorphism of $\pi_{1}(S)$
defined up to an inner automorphism. Moreover, if $\varphi_{1}$ and
$\varphi_{2}$ are isotopic, the corresponding automorphisms of $\pi_{1}(S)$
differ by an inner automorphism. Thus, there is a natural group homomorphism
$\Mod(S)\to\Out(\pi_{1}(S))$. 

The following theorem (due to Dehn, Nielsen \cite{Nielsen1} and Baer
\cite{Baer} for closed surfaces, and to Magnus \cite{Magnus} and
Zieschang \cite{Zieschang} for surfaces with boundary), allows to
obtain the descriptions of the subgroups of outer automorphisms of
the fundamental group corresponding to various mapping class groups
under the map $\Mod(S)\to\Out(\pi_{1}(S))$ described above. 

\begin{thm}
\label{thm:Dehn-Nielsen-Baer}1. If $S$ is a closed connected orientable
surface and is not a sphere, then there is a group isomorphism $\Mod(S)\cong\Out(\pi_{1}(S))$. 

2. If $S$ has non-empty boundary $\partial S=\bigcup_{i=1}^{k}T_{i}\neq\emptyset$
and $\chi(S)<0$ (i.e., $S$ is not a disc or a cylinder), then \begin{eqnarray}
\Mod(S) & \cong & \{\alpha\in\Out(\pi_{1}(S)):\exists\tilde{\alpha}\in\Aut(\pi_{1}(S)),\ [\tilde{\alpha}]=\alpha,\textrm{ }\label{mod-1}\\
 &  & \quad\textrm{and }\forall i\ \ \exists j\ \textrm{s.t. }\tilde{\alpha}([T_{i}])=[T_{j}]\},\label{mod-2}\end{eqnarray}
i.e., the subgroup consisting of outer automorphisms which preserve
the oriented peripheral structure. Moreover,\begin{eqnarray}
\PMod(S) & \cong & \{\alpha\in\Out(\pi_{1}(S)):\forall i\ \ \exists\tilde{\alpha_{i}}\in\Aut(\pi_{1}(S)),\ \textrm{ }\label{pmod-1}\\
 &  & \quad[\tilde{\alpha}_{i}]=\alpha\ \ \ \textrm{s.t. }\tilde{\alpha}([T_{i}])=[T_{i}]\}.\label{pmod-2}\end{eqnarray}

\end{thm}
We will use the relation of mapping class groups with outer automorphisms
of the fundamental group when dealing with the Picard group and its
subgroups in the case of an open symplectic surface.

\begin{rem}
\label{rem:sphere_disc_annulus}Note that
\begin{enumerate}
\item If $S$ is a sphere, then $\Mod(S)\cong\mathbb{Z}_{2}$ and is generated
by an orientation-reversing diffeomorphism. However, $\Out(\pi_{1}(S))=\{ e\}$.
\item If $S$ is a disc, then $\Mod(S)\cong\mathbb{Z}_{2}$, but $\Out(\pi_{1}(S))\cong\{ e\}\cong\PMod(S)$. 
\item If $S=C=I\times S^{1}$ is a cylinder with coordinates $(r,\theta)$,
$r\in I=[1,2]$, then $\PMod(C)\cong\{ e\}$ and $\Out(\pi_{1}(C))\cong\mathbb{Z}_{2}$.
The group $\Mod(C)\cong\mathbb{Z}_{2}\times\mathbb{Z}_{2}$ is generated
by the automorphisms $\Phi_{1}(r,\theta)=(r,-\theta)$ (which induces
the unique non-trivial outer automorphism of $\pi_{1}(C)$) and $\Phi_{2}(r,\theta)=(3-r,\theta)$
(which switches the boundary components of the cylinder). 
\end{enumerate}
\label{rem:M_disconnected-surface}Note that if $S$ is disconnected,
one can still define the mapping class groups in a similar way. In
particular, for a disconnected surface $S$ with connected components
$S_{i}$, $i\in I$ we have\begin{equation}
\mathcal{M}(S)=\Pi_{i\in I}\mathcal{M}(S_{i}).\label{eq:M_script_disconnected_surface}\end{equation}

\end{rem}

\subsection{Fundamental groups of surfaces.}

Recall that an orientable surface of genus $g$ with $b$ boundary
components can be obtained from a $4g$-gon by an identification of
sides according to the word $(a_{1}b_{1}a_{1}^{-1}b_{1}^{-1})\ldots(a_{g}b_{g}a_{g}^{-1}b_{g}^{-1})$
and removal of $b$ disjoint discs from its interior. This leads to
the following presentation of $\pi_{1}(S)$:\begin{eqnarray*}
\pi_{1}(S) & = & \langle a_{1},b_{1},\ldots a_{g},b_{g},\ \  d_{1},\ldots,d_{b}\ :\\
 &  & \  a_{1}b_{1}a_{1}^{-1}b_{1}^{-1}\ldots a_{g}b_{g}a_{g}^{-1}b_{g}^{-1}d_{1}\ldots d_{b}=1\rangle.\ \end{eqnarray*}
Note that if $\partial S\neq\emptyset$ (i.e., $b\neq0$), we can
eliminate a $d_{i}$ from this presentation, and consider $\pi_{1}(S)$
as a free group $\mathbb{F}_{2g+b-1}$ on $(2g+b-1)$ generators.
In particular, we have:

\begin{lem}
\label{lem:commutation}Let $S$ be an orientable surface with a non-trivial
boundary, and let $\gamma\in\pi_{1}(S)$ be the class of a connected
component of $\partial S$. If $\alpha\in\pi_{1}(S)$ commutes with
$\gamma$, then $\alpha=\gamma^{k}$ for some $k\in\mathbb{Z}$.
\end{lem}
\begin{proof}
In our identification of $\pi_{1}(S)$ with the free group $\mathbb{F}_{2g+b-1}$,
the class $\gamma$ corresponds to one of the generators $d_{i}$.
If an element $w\in\mathbb{F}_{2g+b-1}$ commutes with $d_{i}$, it
must be a power of that element.
\end{proof}

\subsection{Dehn twists.}

Observe that if two diffeomorphisms in $\Diff(S\fix\partial S)$ are
isotopic in $\Diff(S\fix\partial S)$ they are also isotopic in $\Diff(S)$.
Thus, there is a natural surjective map $\mathcal{M}(S)\to\PMod(S)$.
The kernel of this map can be described in terms of special elements
in \emph{$\mathcal{M}(S)$}, called the \emph{Dehn twists around the
boundary. }

We will first describe the Dehn twists of an annulus. Let $A=\{(r,\theta):r\in[1,2]\}$
be an annulus with the boundary $\partial A=T_{1}\cup T_{2}$ consisting
of two circles, $T_{1}=\{ r=1\}$ and $T_{2}=\{ r=2\}$. Let $f:\mathbb{R}\to\mathbb{R}$
be a smooth function such that

\begin{itemize}
\item $f(x)=0$ for $x\leq1$;
\item $f(x)=2\pi$ for $x\geq2$;
\item $f'(x)\geq0$;
\end{itemize}
The standard twist automorphism of $A$ is defined by the formula\label{lem:annulus}\[
\Phi(r,\theta)=(r,\theta+f(r)).\]
Note that $\Phi|_{\partial A}=\id$ and the class of $\Phi$ up to
isotopy fixed on the boundary is independent of the choice of $f$.
This class $[\Phi]_{\mathcal{M}(S)}$ is called the (left) \emph{Dehn
twist} of the annulus (see \cite{Dehn-collected-papers} for an English
translation of Dehn's original papers where these twists were introduced).
Moreover, any class in $\mathcal{M}(A)$ is a power of the Dehn twist:

\begin{lem}
The group $\mathcal{M}(A)=\pi_{0}(\Diff(A\fix\partial A))$ is isomorphic
to the infinite cyclic group $\mathbb{Z}$ generated by the Dehn twist. 
\end{lem}
Using orientation-preserving embeddings of the annulus into a surface,
one can transplant the standard twist diffeomorphism of the annulus
to a diffeomorphism of the surface. If $e:A\to S$ is an orientation-preserving
embedding, take the diffeomorphism $e\circ\Phi\circ e^{-1}:\  e(A)\to e(A)$
and extend it by the identity to a diffeomorphism $\Phi_{e}\in\Diff(S)$.
Up to an isotopy fixed on the boundary the diffeomorphism $\Phi_{e}$
depends only on the isotopy class of the embedding, which in turn
is determined by the isotopy class of the oriented image $e(a)$ of
the axis $a=\{ r=3/2\}$ of the annulus. One can call the diffeomorphism
$\Phi_{e}$ a \emph{twist about the circle $T=e(a)$ supported on
the annulus $e(A)$.} The Dehn twist about the circle $T$ is the
class of this diffeomorphism up to an isotopy by diffeomorphisms equal
to the identity outside of the support of the diffeomorphism.

The main property of Dehn twists which we will need in Section \ref{sec:trivializations}
is the following

\begin{lem}
\label{lem:Dehn_twists}Let $S$ be an open orientable surface with
$\partial S=\bigcup_{i=1}^{k}T_{i}$. Let $D(\partial S)\cong\mathbb{Z}^{k}$
be the subgroup of $\mathcal{M}(S)$ generated by the Dehn twists
around the curves parallel to the boundary components. Then there
is a (split) exact sequence\[
D(\partial S)\stackrel{i}{\to}\mathcal{M}(S)\stackrel{p}{\to}\PMod(S),\]
where $i$ is the inclusion, and $p$ is the natural projection, sending
the class of a diffeomorphism fixing the boundary in $\mathcal{M}(S)$
to its class in $\PMod(S)$. 
\end{lem}
\begin{rem}
For $S$ with $\partial S=\cup_{i=1}^{k}T_{i}$, the group $\mathcal{M}(S)$
has two natural subgroups: one is $D(\partial S)\cong\mathbb{Z}^{k}$
generated by the Dehn twists around the curves parallel to the boundary
components, and the other is $G(S)$, generated by the Dehn twists
around the non-separating curves (see, e.g., \cite{MappingClassGroups}).
In fact, $\mathcal{M}(S)$ is generated by these two subgroups. Moreover,
$G(S)$ and $D(\partial S)$ commute, so that $\mathcal{M}(S)\cong G(S)\oplus D(\partial S)\cong\PMod(S)\oplus\mathbb{Z}^{k}$. 
\end{rem}

\subsection{Moser's argument.}

Finally, when we are dealing with a symplectic surface, it is useful
to represent the classes in mapping class groups by symplectomorphisms.
For a surface $S$ with boundary $\partial S=\bigcup_{i=1}^{k}T_{i}$
and a symplectic structure $\omega$, denote by $\Symp(S)$ the group
of symplectomorphisms of $S$, and by $\Symp(S\fix\partial S)$ ---
the subgroup of symplectomorphisms trivial on a neighborhood of the
boundary. Moser's (\cite{Moser}) type argument (extended to noncompact
manifolds by Greene and Shiohama, \cite{Greene-Shiohama}) implies
the following

\begin{lem}
\label{lem:symplectic_repr}Let $(S,\omega)$ be a symplectic surface
with non-trivial boundary, such that any neighborhood of the boundary
has an infinite volume. Let $\alpha\in\Mod(S)$ be a class in the
mapping class group which can be represented by an orientation-preserving
diffeomorphism. Then
\begin{enumerate}
\item There is a symplectomorphism $\varphi\in\Symp(S)$ such that $[\varphi]_{\Mod(S)}=\alpha$.
\item Moreover, if $\alpha\in\PMod(S)$, one can choose $\varphi$ to be
in $\Symp(S\fix\partial S)$. 
\end{enumerate}
\end{lem}
\begin{rem}
In general, $\Mod(S)$ also contains orientation-reversing diffeomorphisms,
which are of course not isotopic to symplectomorphisms. In this case,
the analog of Lemma \ref{lem:symplectic_repr} states that every element
of $\Mod(S)$ is isotopic to a symplectomorphism or an anti-symplectomorphism.
\end{rem}

\section{Symplectic manifolds\label{sec:trivializations}}

\subsection{The Picard group.}

The Picard group of a symplectic manifold was computed by Bursztyn
and Weinstein in \cite{BW1}. For convenience, we outline in this
subsection one of the possible approaches to the computation of this
group.

Let $S$ be a symplectic manifold, $\Gamma(S)$ be its symplectic
groupoid and $X\in\Pic(S)$ be an invertible $S,S$-bimodule. Let
$p\in S$ be a point, and consider the set $X_{p}=s_{X}^{-1}(p)\cap t_{X}^{-1}(p)$,
which we will call the isotropy of $X$ at $p$. Since $X$ is invertible,
the isotropy group $\Gamma_{p}(S)=s^{-1}(p)\cap t^{-1}(p)\cong\pi_{1}(S,p)$
acts freely transitively on $X_{p}$ on the left and on the right.
Thus, $X_{p}$ is a discrete set isomorphic to $\pi_{1}(S)$. Moreover,
for any a fixed $x\in X_{p}$ and $\gamma\in\Gamma_{p}(S)$ there
exists a unique element $\Hol_{x}(\gamma)\in\Gamma_{p}(S)\cong\pi_{1}(S,p)$
such that\begin{equation}
\gamma\cdot x=x\cdot\Hol_{x}(\gamma).\label{eq:def_Holonomy}\end{equation}
Thus, we obtain a map $X_{p}\to\Aut(\pi_{1}(S,p))$ given by $x\mapsto\Hol_{x}\in\Aut(\pi_{1}(S))$.
It turns out that the class of the resulting automorphism $\Hol_{x}$
in the group of outer automorphisms of $\pi_{1}(S)$ is independent
of the choices of $p\in S$ and $x\in X_{p}$. Thus, there is a map\begin{eqnarray*}
\Pic(S) & \to & \Out(\pi_{1}(S)),\\
X & \mapsto & [\Hol_{x}]_{\Out(\pi_{1}(S)).}\end{eqnarray*}
Bursztyn and Weinstein showed in \cite{BW1} that this map is actually
a group isomorphism:

\begin{thm}
\label{thm:Pic_Symplectic}\cite{BW1} For a connected symplectic
manifold $S$, the Picard group is isomorphic to the group of outer
automorphisms of the fundamental group,\begin{equation}
\Pic(S)\cong\Out(\pi_{1}(S)).\label{eq:bw-iso}\end{equation}

\end{thm}
In the remainder of this section, we will describe (certain subgroups
of) the Picard group for symplectic surfaces in terms of mapping class
group.

\subsection{Closed symplectic surfaces.}

In the case that $S$ is a connected closed surface and not a sphere
, the Dehn-Nielsen-Baer's theorem (Theorem \ref{thm:Dehn-Nielsen-Baer})
states that\[
\Out(\pi_{1}(S))\cong\Mod(S).\]
Combining this with (\ref{eq:bw-iso}), we obtain that\[
\Pic(S)\cong\Mod(S),\]
(see also a Remark 6.4 in \cite{BW1}).

Using this isomorphism of the Picard group with a mapping class group,
we can characterize the bimodules in the Picard group which come from
Poisson diffeomorphisms. Let $j:\Poiss(S)\to\Pic(S)$ be the natural
map given by $\varphi\mapsto\Gamma_{\varphi}(S)$. It is not hard
to see that the composition of maps\[
\Poiss(S)\stackrel{j}{\to}\Pic(S)\to\Out(\pi_{1}(S))\to\Mod(S)\]
simply takes a symplectomorphism $\varphi\in\Symp(S)$ to its isotopy
class in $\Mod(S)$. Thus the image of $\Poiss(S)$ in $\Mod(S)$
consists exactly of those elements which can be represented by orientation-preserving
diffeomorphisms. In the isomorphism $\Mod(S)\cong\Out(\pi_{1}(S))$,
these correspond to the subgroup $\Out^{+}(\pi_{1}(S))$ of automorphisms
of $\pi_{1}(S)$ that act trivially on the second cohomology group
$H^{2}(\pi_{1}(S),\mathbb{Z})\cong H^{2}(S,\mathbb{Z})$. Thus we
have the following relations between the group of Poisson diffeomorphisms,
the Picard group and outer automorphisms of the fundamental group
in the case of a closed surface which is not a sphere:

\[
j(\Poiss(S))\cong\Out^{+}(\pi_{1}(S))\subsetneq\Out(\pi_{1}(S))\cong\Mod(S)\cong\Pic(S).\]
(For the sphere, $j(\Poiss(S))\cong\Out(\pi_{1}(S)=\{ e\}\subset\mathbb{Z}_{2}=\Mod(S)\cong\Pic(S)$).
On the other hand, since each class in $\Mod(S)$ for a closed surface
can be represented by either a symplectomorphism or an anti-symplectomorphism,
it follows that every bimodule in the Picard group comes from either
a symplectomorphism or an anti-symplectomorphism.

\subsection{Open surfaces: bimodules trivial near the boundary\label{sub:relPic}}

Let $S$ be a $2$-dimensional open symplectic manifold whose boundary
is non-trivial and consists of a finite number of simple closed curves,
$\partial S=\bigcup_{i=1}^{k}T_{i}$. Let $C_{i}$ be an annular collar
near $T_{i}$, such that $\partial C_{i}=T_{i}\cup T_{i}'$, where
$T_{i}'$ is a curve parallel to $T_{i}$. Let $C=\cup_{i=1}^{k}C_{i}$.
Assume that the symplectic structure is such that each neighborhood
of a boundary curve has an infinite symplectic volume (in particular,
each $C_{i}$ has an infinite volume). Since $\partial S$ is non-trivial,
every diffeomorphism is orientation-preserving, and, thus, by Moser's
argument (see Lemma \ref{lem:symplectic_repr}), is isotopic to a
symplectomorphism. 

In this subsection, we will characterize the image of the map $j:\Poiss(S)\to\Pic(S)$
for the case of an open symplectic surface. Fix once and for all points
$p_{i}\in C_{i}$, $i=1,\ldots,k$. Let \begin{eqnarray*}
\Pic(S,\partial S) & = & \{ X\in\Pic(S)\ :\ \forall i\,\exists x_{i}\in X_{p_{i}}\\
 &  & \quad\textrm{s.t. }\Hol_{x_{i}}([T_{i}])=[T_{i}]\}\end{eqnarray*}
 be the subgroup of the Picard group of $S$ consisting of the bimodules
which are {}``trivial near the boundary'' (for each $i=1,\dots,k$,
the point $x_{i}\in X_{p_{i}}$ can be chosen in such a way that the
induced holonomy automorphism $\Hol_{x_{i}}$ preserves the class
of the corresponding boundary curve). 

\begin{prop}
The subgroup of bimodules in the Picard group of $S$ which come from
Poisson automorphisms is isomorphic to the pure mapping class group,\[
j(\Poiss(S))=\Pic(S,\partial S)\cong\PMod(S).\]

\end{prop}
\begin{proof}
Recall that $\PMod(S)$ is the subgroup of $\Mod(S)$ generated by
the diffeomorphisms which preserve setwise the homotopy classes of
all boundary components of $S$: \begin{eqnarray*}
\PMod(S) & = & \{\alpha\in\Mod(S):\forall T\in\partial S,\exists\tilde{\alpha}\in\Diff(S)),\ [\tilde{\alpha}]=\alpha,\\
 &  & \quad\textrm{s.t. }\tilde{\alpha}(T)=T\}.\end{eqnarray*}
Under the map $\Mod(S)\to\Out(\pi_{1}(S))$, this corresponds to the
following description of the pure mapping class group:\begin{eqnarray*}
\PMod(S) & \cong & \{\alpha\in\Out(\pi_{1}(S))\ :\ \forall T\in\partial S\ \exists\tilde{\alpha}\in\Aut(\pi_{1}(S))\\
 &  & \,\,\,\,\textrm{s.t. $\tilde{\alpha}([T])=[T]\}.$}\end{eqnarray*}
(For surfaces with boundary having negative Euler characteristic this
follows from the Dehn-Nielsen-Baer's theorem, and for the disc and
the cylinder it can be verified directly, see Remark \ref{rem:sphere_disc_annulus}).

Let $w:\Pic(S)\to\Out(\pi_{1}(S))$ be the isomorphism of Theorem
\ref{thm:Pic_Symplectic}, given by\[
w(X)=[\Hol_{x}]_{\Out(\pi_{1}(S))},\quad x\in X\in\Pic(S).\]
We claim that the image of $\Pic(S,\partial S)$ under this map is
isomorphic to $\PMod(S)$.

Let $X\in\Pic(S,\partial S)$ be a bimodule. Then for any $i=1,\ldots,k$,
there exists $x_{i}$ such that $s_{X}(x_{i})=t_{X}(x_{i})=p_{i}\in C_{i}$
and the induced holonomy action is trivial on the homotopy class $[T_{i}]$
of the curve $T_{i}$. Then $\Hol_{x_{i}}\in\Aut(\pi_{1}(S))$ is
such that $[\Hol_{x_{i}}]_{\Out(\pi_{1}(S))}=w(X)$ an$\Hol_{x_{i}}([T_{i}])=[T_{i}].$Hence
$w(X)\in\PMod(S)$.

Conversely, let $\alpha\in\PMod(S)$. Let $\varphi\in\Symp(S,\partial S)$
be a symplectomorphism trivial near the boundary and representing
the class $\alpha$, i.e., $\varphi|_{C_{i}}=\id$, $[\varphi]_{\PMod(S)}=\alpha$.
Let $X$ be the $S,S$-bimodule $\Gamma_{\varphi}(S)$ obtained from
the symplectic groupoid $\Gamma(S)$ by composing the target map with
the symplectomorphism $\varphi$. Thus $s_{X}=s$ and $t_{X}=\varphi\circ t$.
Note that as left $S$-modules (and in particular, as sets) $X=\Gamma(S)$.
Let $\varepsilon:S\to\Gamma(S)$ be the identity bisection of $\Gamma(S)$,
so that $\varepsilon(p)$ is the unit element of $\Gamma(S)$ over
$p\in S$. Let $x_{i}=\varepsilon(p_{i})\in X$, $i=1,\ldots,k$.
Since $\varphi|_{C_{i}}=\id$, $s_{X}(x_{i})=t_{X}(x_{i})=p_{i}$. 

We claim that $\Hol_{x_{i}}([T_{i}])=[T_{i}].$ To see this, let $\tilde{T}_{i}\subset C_{i}$
be a curve parallel to $T_{i}$. Note that $\varepsilon(\tilde{T_{i}})$
is a lift of $\tilde{T_{i}}$, lying in the isotropy of $X$. That
is, $s_{X}(x)=t_{X}(x)$ for all $x\in\varepsilon(\tilde{T}_{i})$.
Since for each $p\in S$, the isotropy $X_{p}=s_{X}^{-1}(p)\cap t_{X}^{-1}(p)$
is discrete, it follows that $\varepsilon(\tilde{T}_{i})$ must be
a horizontal lift of $\tilde{T}_{i}$ for the connection on $X$ that
defines $\Hol_{x_{i}}$. Thus $\Hol_{x_{i}}([T_{i}])=[T_{i}]$ and
$X\in\Pic(S,\partial S)$.

Lastly, we have seen that every bimodule in $\Pic(S,\partial S)$
is in the image of $j:\Poiss(S)\to\Pic(S)$. Conversely, every bimodule
in the image of $j$ must correspond to an element of $\PMod(S)\subset\Out(\pi_{1}(S))$
and thus lie in $\Pic(S,\partial S)$. Thus the image of $j$ is exactly
$\Pic(S,\partial S)\subset\Pic(S)$.
\end{proof}
Thus we have the following relations between the group of Poisson
automorphisms, the Picard group and mapping class groups in the case
of an open surface (which is not a disc or a cylinder):\[
j(\Poiss(S))=\Pic(S,\partial S)\cong\PMod(S)\subsetneq\Mod(S)\subsetneq\Out(\pi_{1}(S))\cong\Pic(S).\]
(If $S$ is a disc, the Picard group is trivial; if $S$ is a cylinder,
$j(\Poiss(S))\cong\PMod(S)=\{ e\}\subset\mathbb{Z}_{2}\cong\Out(\pi_{1}(S))\cong\Pic(S))$.

\subsection{The group of bimodules with chosen trivializations near the boundary.}

For the purposes of computation of the (static) Picard group of a
TSS in Section \ref{sub:Static_Pic_TSSS}, it is useful to consider
the group $\mathcal{P}(S)$ of self-equivalence bimodules over an
open symplectic surface which are trivial near the boundary, together
with a chosen trivialization. By this we mean the following. 

Let $X\in\Pic(S,\partial S)$ be a bimodule trivial near the boundary,
and $x_{i}\in X_{p_{i}},\  i=1,\dots,k$, be the points as in the
definition of $\Pic(S,\partial S)$. One can think of a pair $(X,(x_{1},\dots,x_{k}))$
as a self-equivalence bimodule in $\Pic(S,\partial S)$ with a chosen
trivialization, given by $x_{i}$'s. We say that $(X,(x_{1},\ldots,x_{k}))$
and $(X',(x_{1}',\ldots,x_{k}'))$ are isomorphic if there is a bimodule
isomorphism $f:X\to X'$ preserving the trivialization, $f(x_{j})=x_{j}'$
for all $j=1,\ldots,k$.

\begin{defn}
Let $\mathcal{P}(S)$ be the set of isomorphism classes of pairs $(X,(x_{1},\ldots,x_{k}))$,
where 
\begin{enumerate}
\item $X\in\Pic(S,\partial S)$;
\item $x_{i}\in X_{p_{i}}$, where $p_{i}\in C_{i}$ are fixed points near
the boundary.
\item $\Hol_{x_{i}}([T_{i}])=[T_{i}]$. 
\end{enumerate}
\end{defn}
It is easy to verify the following

\begin{lem}
$\mathcal{P}(S)$ is a group with respect to the following operations:
\begin{enumerate}
\item The multiplication is given by the relative tensor product:\[
(X,(x_{1},\ldots,x_{k}))\times(X',(x_{1}',\ldots,x_{k}'))=(X\otimes_{S}X',\ ([(x_{1},x_{1}')],\ldots,[(x_{n},x_{n}')]),\]
where $[(x_{i},x_{i}')]$ is the equivalence class of $(x_{i},x_{i}')\in X\times_{S}X'$
in $X\otimes_{S}X'$.
\item The inversion is defined by \[
(X,(x_{1},\ldots,x_{k}))^{-1}=(X^{-1},(x_{1},\ldots,x_{k})).\]

\item The identity of $\mathcal{P}(S)$ is the pair $(\Gamma(S),(\varepsilon(p_{1}),\ldots,\varepsilon(p_{k}))$,
where $\varepsilon$ is the identity bisection of the symplectic groupoid
$\Gamma(S)$.
\end{enumerate}
\end{lem}
The main result of this section is that $\mathcal{P}(S)$ is isomorphic
to the group $\mathcal{M}(S)=\pi_{0}(\Diff(S),\partial S)$ of diffeomorphisms
of $S$ fixing a neighborhood of the boundary up to isotopies:

\begin{prop}
\label{pro:PandM}$\mathcal{P}(S)\cong\mathcal{M}(S)$. 
\end{prop}
\begin{proof}
Let $\alpha\in\mathcal{M}(S)$ be an isotopy class and $\varphi\in\Symp(S,\partial S)$
be a symplectomorphism trivial near the boundary which represents
this class, $\varphi|_{C}=\id$, $[\varphi]_{\mathcal{M}(S)}=\alpha$.
Consider the map $\theta:\mathcal{M}(S)\to\mathcal{P}(S)$ given by\begin{equation}
[\varphi]\mapsto(\Gamma_{\varphi}(S),\ (\varepsilon(p_{1}),\ldots,\varepsilon(p_{n})).\label{eq:map_theta}\end{equation}

\begin{claim}
The map $\theta$ is well-defined.
\end{claim}
\begin{proof}
Let $\alpha\in\mathcal{M}(S)$. We need to show that $\theta(\alpha)$
does not depend on the choice of $\varphi$ such that $[\varphi]=\alpha$.
Since $\theta$ clearly takes compositions of symplectomorphisms to
products of the corresponding elements in $\mathcal{P}(S)$, it is
sufficient to prove this for a representative of the class of the
identity diffeomorphism.

Let $\varphi\in\Symp(S,\partial S)$ be a symplectomorphism such that
$[\varphi]=\id\in\mathcal{M}(S)$. Thus there exists an isotopy $\varphi^{(\tau)},\,\tau\in[0,1],$
between the identity diffeomorphism $\id=\varphi^{(0)}$ and $\varphi=\varphi^{(1)}$.
We may assume that $\varphi^{(\tau)}\in\Symp(S,\partial S)$ for all
$\tau\in[0,1]$. Let $X^{(\tau)}$ be the bimodule $X^{(\tau)}=\Gamma_{\varphi^{(\tau)}}(S)$. 

Let $x\in X^{(0)}=\Gamma(S)$ be the homotopy class of a path $\gamma$,
so that $s(x)=\gamma(0)$ and $t(x)=\gamma(1)$. Let $y^{(\tau_{0})}\in X^{(\tau_{0})}=\Gamma_{\varphi^{(\tau_{0})}}(S)$
be the homotopy class of a path $\beta^{(\tau_{0})}:[0,\tau_{0}]\to S$
given by\begin{equation}
\beta^{(\tau_{0})}(\tau)=\left(\varphi^{(\tau)}\right)^{-1}(\gamma(1)),\label{eq:beta_tau_0}\end{equation}
so that\begin{eqnarray}
 &  & s_{X^{(\tau_{0})}}(y^{(\tau_{0})})=s(y^{(\tau_{0})})=\gamma(1),\label{eq:s_tau}\\
 &  & t_{X^{(\tau_{0})}}(y^{(\tau_{0})})=\varphi^{(\tau_{0})}(t(y^{(\tau_{0})}))=\gamma(1).\label{eq:t_tau}\end{eqnarray}
For each $\tau_{0}\in[0,1]$, let $\alpha^{(\tau_{0})}:X^{(0)}\to X^{(\tau_{0})}$
be the map $x\mapsto y^{(\tau_{0})}\circ x$. Then\begin{eqnarray*}
s_{X^{(\tau_{0})}}(\alpha^{(\tau_{0})}(x)) & = & s(x),\\
t_{X^{(\tau_{0})}}(\alpha^{(\tau_{0})}(x)) & = & t(x).\end{eqnarray*}
 This implies that $\alpha^{(\tau_{0})}$ sends the identity bisection
of the symplectic groupoid $X^{(0)}$ to the identity bisection of
the bimodule $X^{(\tau_{0})}$, and thus defines a bimodule isomorphism
$X^{(0)}\cong X^{(\tau_{0})}$. Since $\varphi^{(\tau)}|_{C_{i}}=\id|_{C_{i}}$,
it follows that this isomorphism preserves the chosen trivialization,
$\alpha^{(\tau_{0})}(x_{i})=x_{i}$, $\forall i=1,\ldots,k$. Thus
the map $\theta$ is indeed well-defined.
\end{proof}
Let $q:\mathcal{P}(S)\to\Pic(S,\partial S)$ be the projection map
taking the pair $(X,(x_{1},\ldots,x_{k}))$ to $X$, and let $\pr:\mathcal{M}(S)\to\PMod(S)$
be the natural map that takes the class of a diffeomorphism $\varphi$
in $\mathcal{M}(S)$ to its class in $\PMod(S)$. Note that both maps
are surjective group homomorphisms. We then have the following commutative
diagram:\begin{equation*} \begin{CD} \mathcal{M}(S) @>\theta>>\mathcal{P}(S) \\ @V{\pr}VV @V{q}VV \\ \PMod(S) @>w^{-1}>> \Pic(S,\partial S). \end{CD}\end{equation*}
\begin{claim}
$\theta$ is onto. 
\end{claim}
\begin{proof}
Since $\pr$ and $q$ are onto and $w$ is an isomorphism, it is sufficient
to prove that $\textrm{ker}(q)\subseteq\textrm{Im}(\theta)$.

Let $(X,(x_{1},\ldots,x_{k}))\in\ker q$. Thus we may assume that
$X$ is the identity bimodule $\Gamma(S)$ and that $x_{1},\ldots,x_{k}\in\Gamma(S)$
are points satisfying $s(x_{j})=t(x_{j})=p_{j}$, $j=1,\ldots,k$.
Let $\gamma_{j}$ be a curve based at $p_{j}$ whose homotopy class
is $x_{j}$. The condition that\[
\Hol_{x_{j}}([T_{j}])=([T_{j}])\]
means by definition that if $\lambda_{j}$ is a curve based at $p_{j}$
and parallel to $T_{j}$, then $\gamma_{j}\lambda_{j}\gamma_{j}^{-1}$
has the same class in $\pi_{1}(S)$ as $\lambda_{j}$. By Lemma \ref{lem:commutation},
the only elements of the fundamental group of $S$ that commute with
the class of a boundary curve $[\lambda_{j}]=[T_{j}]$ are the powers
of the class of this boundary curve. It follows that $x_{j}$ must
be represented by a power of $\lambda_{j}$, i.e., $x_{j}=\left[\lambda_{j}\right]^{n_{j}}$
for some $n_{j}\in\mathbb{Z}$. Thus $(X,(x_{1},\ldots,x_{n}))\cong(\Gamma(S),([T_{1}]^{n_{1}},\ldots,[T_{k}]^{n_{k}}))$
for some integers $n_{1},\ldots,n_{k}$.

We will now construct a class in $\mathcal{M}(S)$ which is mapped
to this bimodule by $\theta$. Let $D_{i}$ be a cylindrical neighborhood
parallel to $C_{i}$ and adjacent to $C_{i}$, and let $S'=S\setminus(\bigcup_{i}C_{i}\cup\bigcup_{i}D_{i})$.
Let $\varphi:S\to S$ be a symplectomorphism with the following properties:
\begin{itemize}
\item $\varphi|_{D_{i}}$ is (homotopic to) the $n_{j}$-th power of the
Dehn twist on $D_{i}$ , $\forall i=1,\ldots,k$;
\item $\varphi$ is trivial outside of $\bigcup D_{i}$, i.e., $\varphi|_{S\setminus\bigcup_{i}D_{i}}=\id$.
\end{itemize}
We claim that $\theta([\varphi])=(X,(x_{1},\ldots,x_{k}))$. Indeed,
let $\varphi^{(\tau)},\tau\in[0,1]$ be an isotopy between $\id=\varphi^{(0)}$
and $\varphi=\varphi^{(1)}$ . Of course $\varphi^{(\tau)}$ no longer
preserve pointwise the boundary of $S$ for $0<\tau<1$; in fact,
for $\tau\in[0,1]$, the map $\tau\mapsto\varphi^{(\tau)}(p_{i})$
traces out a curve which is homotopic to $T_{i}^{n_{i}}$.

For each $\tau\in[0,1]$, let $X^{(\tau)}=\Gamma_{\varphi^{(\tau)}}(S)$
and let $\alpha^{(\tau_{0})}:X^{(0)}\to X^{(\tau_{0})}$ be the bimodule
isomorphism described in the proof of Claim 1. Identifying $\Gamma(S)$
with the space of homotopy classes of paths in $X$, we view $\alpha^{(\tau_{0})}(x_{j})$
as the class of the path $\tau\mapsto\varphi^{(\tau)}(p_{j})$, $0\leq\tau\leq\tau_{0}$.
It follows that $\alpha^{(1)}(x_{j})$ is homotopic to $T_{j}{}^{n_{j}}$.
Thus, $\alpha^{(1)}$ is an isomorphism between $(X,(x_{1},\ldots,x_{n}))$
and $\theta([\varphi])=(\Gamma(S),\varepsilon(p_{1}),\ldots,\varepsilon(p_{k}))$.
Hence $\theta$ is onto.
\end{proof}
\begin{claim}
$\theta$ is injective. 
\end{claim}
\begin{proof}
Let $\varphi\in\Symp(S,\partial S)$ be such that $[\varphi]_{\mathcal{M}(S)}\in\ker\theta$.
Since $w$ is an isomorphism, it follows that $[\varphi]_{\mathcal{M}(S)}\in\ker(\pr)$.
We may therefore assume (Lemma \ref{lem:Dehn_twists}) that there
are cylindrical neighborhoods $D_{1},\ldots,D_{k}$ adjacent to $C_{1},\ldots,C_{k}$
and parallel to them, so that $\varphi$ is identity outside of $D_{1}\cup\cdots\cup D_{k}$
and the restriction of $\varphi$ to each $D_{j}$ is the $n_{j}$-th
power of the Dehn twist on $D_{j}$.

We saw earlier that $\theta([\varphi])$ is isomorphic to $(\Gamma(S),([T_{1}]^{n_{1}},\ldots,[T_{k}]^{n_{k}}))$.
Since $[\varphi]\in\ker\theta$, it follows that there is a bimodule
isomorphism $\Phi:\Gamma(S)\to\Gamma(S)$ such that $\Phi([T_{i}]^{n_{i}})=\varepsilon(p_{i})$,
$\forall i=1,\ldots,k$, where $\varepsilon$ denotes the identity
bisection. Since $\varepsilon:S\to\Gamma(S)$ is an identity bisection
of $\Gamma(S)$, $\sigma=\Phi^{-1}\circ\varepsilon:S\to\Gamma(S)$
is also an identity bisection. Thus there is an identity bisection
$\sigma:S\to\Gamma(S)$ passing through $x_{j}=[T_{j}]^{n_{j}}$ over
$p_{j}$ for all $j=1,\ldots,k$. By property 1 of an identity bisection
(Definition \ref{def:identity_bisection}), for any $p\in S$ the
element $\sigma(p)$ lies in the discrete set $s^{-1}(p)\cap t^{-1}(p)=\pi_{1}(S,p)$,
and, therefore, $\sigma(p)=\gamma(p)\cdot\varepsilon(p)$ for some
(locally constant) map $\gamma:S\to\pi_{1}(S)$. By property 2 in
the same definition, the left and right groupoid actions by an element
$\beta\in\Gamma_{p}(S)=\pi_{1}(S,p)$ on $\sigma(p)=\gamma(p)\cdot\varepsilon(p)$
commute, which implies that $\gamma(p)$ lies in the center of $\pi_{1}(S,p)$
for all $p\in S$. Moreover, the map $p\mapsto\gamma(p)$ is globally
constant.

Since $\partial S\neq\emptyset$, the fundamental group of $S$ can
be identified with a free group. Hence $S$ is either a cylinder,
or its fundamental group is a free group on two or more generators
and thus has a trivial center. In the latter case, we get that $\sigma=\varepsilon$
so that $x_{j}=\varepsilon(p_{j})$, $j=1,\ldots,k$, and thus $n_{j}=0$
for all $j$. This implies $[\varphi]=\id\in\mathcal{M}(S)$.

If $S$ is a cylinder, $\partial S=T_{1}\cup T_{2}$. Let $T$ be
a separating circle on the cylinder, so that $[T]$ is a generator
of the fundamental group of $S$. The identity bisections $\varepsilon_{k}:S\to\Gamma(S)$
of the symplectic groupoid of $S$ have the form $\varepsilon_{k}(p)=[T(p)]^{k}$,
where $k\in\mathbb{Z}$, and $T(p)$ is a closed curve based at $p$
and homotopic to $T$. Thus if there is an identity bisection of $X$
through $x_{1}$ and $x_{2}$, it must be that $x_{1}$ and $x_{2}$
are the homotopy classes of $T^{k}$ for the same $k$. In this case
the two Dehn twists making up $\varphi$ cancel and $[\varphi]=\id\in\mathcal{M}(S)$.

This concludes the proof that $\theta$ is injective.
\end{proof}
This completes the proof of the Proposition.
\end{proof}

\subsection{A remark on bisections of bimodules.}

Let $X$ be an $(S,S)$-bimodule.

\begin{defn}
\label{def:bisections}A map $\sigma:S\to X$ is called a \emph{bisection}
of $X$, if
\end{defn}
\begin{enumerate}
\item $s_{X}\circ\sigma=\id$, and $t_{X}\circ\sigma=\varphi$ is a symplectomorphism
of $S$.
\item $\sigma(S)$ is a lagrangian submanifold of $X$.
\end{enumerate}
Note that the first condition means that $\sigma$ is a section of
the source, map, while $\sigma\circ\varphi^{-1}$ is a section of
the target map. 

Let $\Bis(X)$ be the set of all bisections of $X$, and let\[
\mathcal{B}(S)=\{(X,\sigma):X\in\Pic(S),\ \sigma\in\Bis(X)\}\]
be the set of self-equivalence bimodules with chosen bisections. Then
$\mathcal{B}(S)$ is a group under the relative tensor product operation.
The inverse of a pair $(X,\sigma)$ is given by $(X,\sigma=(X^{\textrm{-1}},\sigma\circ\varphi^{-1})$,
where $X^{-1}=X^{\textrm{op}}\ $ is the opposite bimodule, and $\varphi$
is the symplectomorphism given by $t\circ\sigma$, so that $s_{X^{-1}}^{\textrm{}}\circ\sigma=\id$
and $t_{X^{-1}}\circ\sigma=\varphi^{-1}$.

The map $\Phi:\Symp(S)\to\mathcal{B}(S)$ given by $\varphi\mapsto(\Gamma_{\varphi}(S),\varepsilon)$
is a group isomorphism. 

Endow the space of all symplectomorphisms $\Symp(S)$ of $S$ with
the $C^{\infty}$ topology. The space $\mathcal{B}(S)$ can be endowed
with a topology, making this isomorphism a homeomorphism. We will
say that $(X,\sigma)$ and $(X',\sigma')$ are \emph{isotopic} if
there is a continuous path $(X^{(t)},\sigma^{(t)})$ in $\mathcal{B}(S)$,
joining $(X,\sigma)$ and $(X',\sigma')$.

From the definition of $\mathcal{M}(S)$ and the isomorphism between
$\mathcal{B}(S)$ and the group of all symplectomorphisms of $S$
we deduce that\[
\mathcal{B}(S)/\textrm{isotopy}\cong\mathcal{M}(S).\]
Thus our result that $\mathcal{M}(S)\cong\mathcal{P}(S)$ implies
that\[
\mathcal{P}(S)\cong\mathcal{B}(S)/\textrm{isotopy}=\pi_{0}(\Symp(S)).\]

\section{\label{sec:TSS} The Picard Group of a TSS.}

\subsection{Topologically stable structures on compact oriented surfaces.}

From now on, let $P$ be a compact connected oriented surface and
let $\pi$ be a Poisson structure on $P$ with at most linear degeneracies.
The zero set $Z\subset P$ of such a structure consists of a finite
number of simple closed curves, $Z=\cup_{i=1}^{n}T_{i}$. Such structures
are called \emph{topologically stable} structures (or, TSS) since
the topology of their zero set does not change under small perturbation
of the Poisson tensor.

By an easy application of the integrability criteria of Crainic and
Fernandes (\cite{CF1}), we know that TSS are integrable. (Alternatively,
integrability follows from a result of Debord (\cite{Debord}), since
the anchor of the Lie algebroid associated to a TSS is injective on
a dense open set $P\setminus Z$.) 

By the main result of \cite{Radko-classification}, the following
invariants:

\begin{itemize}
\item topological class of the oriented zero set;
\item modular periods around the zero curves;
\item a generalized Liouville volume; 
\end{itemize}
completely classify TSS up to orientation-preserving Poisson diffeomorphisms.
Moreover, according to \cite{BW1} (see also \cite{BR}), the topology
of the oriented zero set together with the modular periods around
the zero curves classify TSS up to Morita equivalence.

\subsection{The Picard group of $(\mathbb{R}^{2},\pi=x\partial_{x}\wedge\partial_{y})$.}

Since a TSS vanishes linearly on a zero curve $T$, in a neighborhood
of a point $p\in T$ it is isomorphic to $\mathbb{R}^{2}$ with the
Poisson structure $\pi=x\partial_{x}\wedge\partial_{y}$. Let us start
with the description of the identity bimodule $\Gamma(\mathbb{R}^{2})$
(i.e., the symplectic groupoid) of this structure. Since $\pi=x\partial_{x}\wedge\partial_{y}$
is the Lie-Poisson structure on the dual of the Lie algebra of the
group $G$ of affine transformation of the line, it follows that $\Gamma(\mathbb{R}^{2})\cong T^{*}G$.
As a manifold, $\Gamma(\mathbb{R}^{2})$ is diffeomorphic to $\mathbb{R}^{4}$
with coordinates $(x,y,p,q)$, in which the Lie groupoid structure
is given by\begin{eqnarray*}
 &  & s((x,y,p,q))=(x,y),\qquad t((x,y,p,q))=(xe^{p},y+xq),\\
 &  & (x,y,p,q)\cdot(x',y',p',q')=(x,y,p+p',q+e^{p}q'),\\
 &  & \textrm{where }x'=xe^{p},\, y'=y+xq,\\
 &  & (x,y,p,q)^{-1}=(xe^{p},y+xq,-p,-qe^{-p}).\end{eqnarray*}
The symplectic form on $\Gamma(\mathbb{R}^{2})$ is given by \begin{eqnarray}
\Omega & = & t^{*}(d(\ln x)\wedge dy)-s^{*}(d(\ln x)\wedge dy))\label{eq:Omega}\\
 & = & -qdx\wedge dp+dx\wedge dq-dy\wedge dp+xdp\wedge dq,\nonumber \end{eqnarray}
and the corresponding Poisson tensor is\begin{equation}
\Pi=-x\partial_{x}\wedge\partial_{y}+\partial_{x}\wedge\partial_{q}-\partial_{y}\wedge\partial_{p}-q\partial_{y}\wedge\partial_{q}.\label{Pi}\end{equation}
Using this description of the symplectic groupoid, we can compute
the Picard group: 

\begin{prop}
\label{pro:Picard_plane}Let $\pi=x\partial_{x}\wedge\partial_{y}$
be a Poisson structure on $\mathbb{R}^{2}$. Then
\begin{enumerate}
\item The static Picard group $\StatPic(\mathbb{R}^{2},\pi)$ is trivial.
\item The full Picard group is given by\[
\Pic(\mathbb{R}^{2},\pi)\cong\mathbb{Z}_{2}\times\mathbb{R}\cong\Out(\mathbb{R}^{2},\pi),\]
where the generator of $\mathbb{Z}_{2}$ corresponds to the flip $(x,y)\mapsto(-x,y)$,
and $t\in\mathbb{R}$ corresponds to a shift $(x,y)\mapsto(x,y+t)$.
\end{enumerate}
\end{prop}
\begin{proof}
Let $(\mathbb{R}^{2},\pi)\stackrel{s_{X}}{\leftarrow}(X,\,\Omega_{X})\stackrel{t_{X}}{\rightarrow}(\mathbb{R}^{2},\pi)$
be a Morita self-equivalence bimodule of $(\mathbb{R}^{2},\pi)$,
inducing the identity map on the leaf space. Define on $X$ the functions
$x_{i}(s)=x(J_{i}(s))$, $y_{i}(s)=y(J_{i}(s))$, $i=1,2$. Let $T=\{(0,y)|\, y\in\mathbb{R}\}$
be the zero set of the Poisson structure. Define\begin{eqnarray}
 &  & p(s)=\textrm{ln}\frac{x_{2}}{x_{1}},\label{p-def}\\
 &  & q(s)=\frac{y_{2}-y_{1}}{x_{1}}.\label{q-def}\end{eqnarray}
 The functions $p$ and $q$ are well-defined on $X\setminus{s_{X}}^{-1}(T)$,
where both $x_{1}$ and $x_{2}$ are nonzero. We claim that $p$ and
$q$ extend smoothly to all of $X$.

To see this for $p$, note that $p=\textrm{ln}x_{2}-\textrm{ln}x_{1}$
is a function whose Hamiltonian vector field $H_{p}=\tilde{\pi}(dp)$
projects by the moment maps to the modular vector field of the Poisson
structure with respect to the standard area form $dx\wedge dy$ on
$\mathbb{R}^{2}$, that is, $(s_{X})_{*}(H_{p})=(t_{X})_{*}(H_{p})=\partial_{y}$.
By \cite{Ginzburg-Lu}, there exists a \emph{smooth} function $h\in C^{\infty}(X)$
which has the same property, $(s_{X})_{*}(H_{h})=(t_{X})_{*}(H_{h})=\partial_{y}$
. Since on the dense subset $X\setminus{s_{X}}^{-1}(T)\subset X$
the map $x\mapsto(s_{X}(x)\ ,t_{X}(x))$ is one-to-one and onto, it
follows that such a vector field is unique, i.e., $H_{p}=H_{h}$ on
$X\setminus{s_{X}}^{-1}(T)$. Thus, $p-h$ is a locally constant function
on $X\setminus s_{X}^{-1}(T)$. Since $x_{2}=x_{1}e^{p}$ is smooth
on $X$, it follows that $p=h+\textrm{const}$, and is therefore,
smooth.

To prove that $q$ extends to all of $X$, it is enough to recall
that $y_{1}=y_{2}$ on $s_{X}^{-1}(T)$, which follows from the assumption
that the map induced by $X$ on the leaf space is the identity. 

Note that because $X\setminus s_{X}^{-1}(T)$ is symplectomorphic
to $(\mathbb{R}^{2}\setminus T)\times(\mathbb{R}^{2}\setminus T)^{\op}$,
it follows that the symplectic form on $X\setminus s_{X}^{-1}(T)$
is given by\begin{equation}
\Omega=-qdx_{1}\wedge dp+dx_{1}\wedge dq-dy_{1}\wedge dp+x_{1}dp\wedge dq,\label{Omega_on_S}\end{equation}
and thus on all of $X$ by continuity.

We claim that $(x_{1},y_{1},p,q)$ is a coordinate system on $X$.
First, since the map $\psi:X\to\Gamma(\mathbb{R}^{2})$ given by $s\mapsto(x_{1}(s),y_{1}(s),p(s),q(s))$
preserves the symplectic form, and, hence, the volume form, it is
a local diffeomorphism. Since $X$ is a Morita self-equivalence bimodule,
there is a diffeomorphism $\Phi:\  X\otimes_{\mathbb{R}^{2}}X^{-1}\to\Gamma(\mathbb{R}^{2})$
. This implies that the map $\psi$ is one-to-one, and, therefore
$(x_{1},y_{1},p,q)$ is a coordinate system on $X$, which establishes
the diffeomorphism of $X$ and $\Gamma(\mathbb{R}^{2})$. Thus, every
invertible bimodule preserving the leaf space pointwise is isomorphic
to the identity bimodule, and, therefore, the static Picard group
is trivial.

To compute the full Picard group, we apply the exact sequence (\ref{eq:Picard_exact_sequence}).
Any automorphism of the leaf space $L(\mathbb{R}^{2})$ that comes
from a Morita self-equivalence bimodule must preserve the modular
vector field (see \cite{Ginzburg-Lu}). Thus, the restriction of an
automorphism of $L(\mathbb{R}^{2})$ to the zero set must be a translation
by the flow of the restriction of a modular vector field. Thus the
image of the map $h:\textrm{Pic}(\mathbb{R}^{2},\pi)\to\textrm{Aut}(L(\mathbb{R}^{2}))\cong\mathbb{Z}_{2}\times\Diff(\mathbb{R})$
is contained in $\mathbb{Z}_{2}\times\mathbb{R}$, where $\mathbb{Z}_{2}$
is generated by the interchange of the two-dimensional leaves, and
$\mathbb{R}$ is generated by a shift along the line of zero-dimensional
leaves. Since any automorphism $(\sigma,t)\in\mathbb{Z}_{2}\times\mathbb{R}$
can be realized by the Poisson automorphism $\theta_{(\sigma,t)}(x,y)=(\sigma x,\  y+t),$
which gives rise to a self-equivalence bimodule $\Gamma_{\theta_{(\sigma,t)}}(\mathbb{R}^{2})$,
it follows that $\textrm{Im}(h)\simeq\mathbb{Z}_{2}\times\mathbb{R}^{2}$.
Since an automorphism $\theta_{(\sigma,t)}$ is inner iff $(\sigma,t)=(1,0)=e_{\mathbb{Z}_{2}\times\mathbb{R}}$,
we conclude that\[
\textrm{Pic}(\mathbb{R}^{2},\pi)\simeq\textrm{Out}(\mathbb{R}^{2},\pi)\simeq\mathbb{Z}_{2}\times\mathbb{R}.\]
 
\end{proof}

\subsection{The Static Picard group\label{sub:Static_Pic_TSSS} of a TSS}

Let $X\in\StatPic(P)$. Let $S\subset P$ be a $2$-dimensional symplectic
leaf. The $X_{|S}\in\Pic(S,\partial S)$, and for $p\in S$ the isotropy
$X_{p}=s_{X}^{-1}(p)\cap t_{X}^{-1}(p)$ is a discrete set, isomorphic
to the fundamental group $\pi_{1}(S,p)$. 

The following lemma shows that locally, in a neighborhood of a zero
curve $T$ of a TSS, there is at most one lift of a curve which crosses
$T$ to the {}``isotropy subbundle'' $\cup_{p\in P}X_{p}$ of the
bimodule. 

\begin{lem}
\label{lemma:repelling}Let $S$ be a $2$-dimensional symplectic
leaf of $P$. Let $\gamma:[0,1]\to P$ be a curve in $P$ such that 
\begin{itemize}
\item $\gamma(t)\in S$ for $t\in[0,1)$;
\item $\gamma(1)\in T\subset\partial S\subset Z.$
\end{itemize}
For $\tau\in[0,1)$, let $\gamma_{\tau}:\ [0,\tau]\to S$ be the curve
$\gamma_{\tau}(t)=\gamma(t)$.

Let $p=\gamma(0)$ and $x_{1},x_{2}\in X_{p}$. For $\tau\in[0,1)$,
define\[
x_{j}(\tau)=[\gamma_{\tau}]\cdot x_{j}\cdot[\gamma_{\tau}]^{-1},\quad j=1,2,\]
where $[\gamma_{\tau}]\in\Gamma(P)$ is the class of the cotangent
path $\tilde{\pi}^{-1}(\gamma_{\tau})$, and $\cdot$ denotes the
left and right actions of $\Gamma(P)$ on $X$.

If both limits $\lim_{\tau\to1}x_{1}(\tau)$ and $\lim_{\tau\to1}x_{2}(\tau)$
exist, then $x_{1}=x_{2}$. 
\end{lem}
\begin{proof}
Let us assume that both limits, $\lim_{\tau\to1}x_{1}(\tau)$ and
$\lim_{\tau\to1}x_{2}(\tau)$, exist.

Since $X$ is invertible, the relative tensor product $X\otimes_{P}X^{\op}$
is isomorphic to the identity bimodule $\Gamma(P)$. Let \begin{eqnarray*}
y_{1}(\tau) & = & (x_{1}(\tau),x_{1}(\tau))\in X\otimes_{P}X^{\op}\cong\Gamma(P),\\
y_{2}(\tau) & = & (x_{2}(\tau),x_{1}(\tau))\in X\otimes_{P}X^{\op}\cong\Gamma(P).\end{eqnarray*}
Then $x_{1}(0)=x_{2}(0)$ if and only if $y_{1}(0)=y_{2}(0)$. Moreover,
$y_{j}(\tau)$ satisfy\begin{equation}
s(y_{j}(\tau))=t(y_{j}(\tau))=\gamma(\tau),\quad j=1,2.\label{eq:s_and_t_y_i}\end{equation}
Let $z(\tau)=y_{1}(\tau)\cdot y_{2}(\tau)^{-1}\in\Gamma(P).$ As usual,
we will view the elements of $\Gamma(P)$ as classes of cotangent
paths up to cotangent homotopy. For each $\tau\in[0,1)$, let $(a^{(\tau)},\eta^{(\tau)}$)
be a cotangent path representing $z(\tau)$. Here $a^{(\tau)}:[0,1]\to T^{*}P$
is a path in the cotangent bundle, $\eta^{(\tau)}:[0,1]\to P$ is
the base path, and the compatibility condition states\begin{equation}
\tilde{\pi}(a^{(\tau)}(t))=\frac{d\eta}{dt}^{(\tau)}(t).\label{eq:cotangent_paths}\end{equation}
 In particular, since $y_{j}(\tau)$ is in the isotropy group of $\Gamma(P)$
at $\gamma(\tau)$, the homotopy class of $\eta^{(\tau)}$ is trivial
if and only if $z(\tau)$ is the identity element, i.e., if and only
if $x_{1}(\tau)=x_{2}(\tau)$.

Endow $P$ with a fixed metric. If the homotopy class of $\eta^{(\tau)}$
is non-trivial, we have\[
\inf_{\tau}\ \ \ \inf_{((a^{(\tau)},\ \eta^{(\tau)})\ \textrm{representing }z(\tau))}\ \ \sup_{t}\Big\Vert\frac{d\eta^{(\tau)}}{dt}(t)\Big\Vert=C>0\]
for the norm $\Vert\cdot\Vert$ on $TP$ coming from our choice of
a metric on $P$. But since $\gamma(1)\in T$ , where the Poisson
structure vanishes, condition (\ref{eq:cotangent_paths}) implies
that\[
\lim_{\tau\to1}\ \ \ \inf_{((a^{(\tau)},\eta^{(\tau)})\textrm{ representing}\  z(\tau))}\ \ \ \sup_{t}\Vert a^{(\tau)}(t)\Vert\to\infty.\]
This contradicts the existence of the limit\[
\lim_{\tau\to1}z(\tau)=\lim_{\tau\to1}y_{1}(\tau)\cdot\lim_{\tau\to1}y_{2}(\tau)^{-1},\]
which follows from the assumption that $\lim_{\tau\to1}x_{i}(\tau)$
exists for $i=1,2$. Thus we must have that $x_{1}=x_{2}$.
\end{proof}
The following corollary is immediate:

\begin{cor}
Let $S$ be a $2$-dimensional symplectic leaf of $P$. Let $\gamma:[0,1]\to P$
be a curve in $P$ such that 
\begin{itemize}
\item $\gamma(t)\in S$ for $t\in[0,1)$;
\item $\gamma(1)\in T\subset\partial S\subset Z.$
\end{itemize}
If $\sigma_{i}:[0,1]\to X$, $i=1,2$, are any two curves in the isotropy
subbundle lying over $\gamma$, i.e., \[
s_{X}(\sigma_{i}(t))=t_{X}(\sigma_{i}(t))=\gamma(t),\quad i=1,2,\]
then $\sigma_{1}=\sigma_{2}$.
\end{cor}
Next, using this corollary, we will show that the restriction of a
static Morita self-equivalence bimodule $X\in\StatPic(P)$ to a neighborhood
of each zero curve has a unique identity bisection. Thus, we will
obtain the local identity bisection of $X$ on a neighborhood of the
zero set.

\begin{lem}
\label{lem:local-bisections-exist}Let $X\in\StatPic(P)$ and $T\subset P$
be a zero curve. Then there is a neighborhood $N$ of $T$ and a map
$\sigma:N\to X$ so that $t_{X}(\sigma(x))=s_{X}(\sigma(x))$ for
all $x\in N$. Moreover, $\sigma$ is uniquely determined by this
property.
\end{lem}
\begin{proof}
First, we will show that there is a unique identity bisection on a
neighborhood of a point on the zero curve. 

Let $p_{0}\in T$ be a point. Since $s_{X}$ is a submersion, we can
find a cross-section defined on a neighborhood $N_{0}'$ of $p_{0}$,
i.e., there is a map $\sigma:N_{0}'\to X$, so that $s_{X}\circ\sigma=\id$.
Since the source and target maps coincide on $X_{|T}$, it follows
that on a neighborhood $N_{0}\subset N_{0}'$, the composition $t_{X}\circ\sigma$
is a diffeomorphism. We can also assume that the Poisson manifold
$(N_{0},\pi|_{N_{0}})$ is isomorphic to $(\mathbb{R}^{2},\  x\partial_{x}\wedge\partial_{y})$. 

To construct an identity bisection of $X_{|N_{0}}$, we would like
to apply the result of Proposition \ref{pro:Picard_plane}, which
states that every bimodule over $(\mathbb{R}^{2},\  x\partial_{x}\wedge\partial_{y})$
is trivial, and therefore, by Lemma \ref{lem:identityBisectionTrivial},
has a unique identity bisection. The problem is that $X_{|N_{0}}$
has disconnected fibers and thus is not a bimodule in our sense. To
find a bimodule $X_{0}\subset X_{|N_{0}}$, identify the symplectic
groupoid $\Gamma(N_{0})$ with the connected component of the identity
bisection in $\Gamma(P)_{|N_{0}}$. Denote by $\cdot$ the left and
right actions of $\Gamma(N_{0})$ on $X$ obtained by restricting
the actions of $\Gamma(P)$. Let\[
X_{0}=\Gamma(N_{0})\cdot\sigma(N_{0})\cdot\Gamma(N_{0}),\]
where $\sigma:N_{0}\to X$ is as above. Then $X_{0}$ is the connected
component of $\sigma(N_{0})$ inside of $X_{|N_{0}}$, and thus a
symplectic manifold of the same dimension. Moreover, $X_{0}$ is clearly
a $\Gamma(N_{0})$-bimodule. Furthermore, since the isotropy groups
of $\Gamma(N_{0})$ are trivial (they are isomorphic to the trivial
fundamental group $\pi_{1}(N_{0})$), it follows that the action of
$\Gamma(N_{0})$ on the fibers of $X_{0}$ is free; it is by definition
transitive. Hence $X_{0}\in\StatPic(N_{0})$. By Proposition \ref{pro:Picard_plane},
$\StatPic(N_{0})=\StatPic(\mathbb{R}^{2})$ is trivial, and thus $X_{0}$
is isomorphic to $\Gamma(\mathbb{R}^{2})$. Identify $X_{0}$ from
now on with $\Gamma(\mathbb{R}^{2})$. Let $\varepsilon_{0}:N_{0}\to\Gamma(N_{0})=X_{0}\subset X$
be the identity bisection, which exists by Lemma \ref{lem:identityBisectionTrivial}.

We will now extend this identity bisection to a cylindrical neighborhood
of $T$. Let $\xi$ be a modular vector field of $\pi$ with the property
that its orbits are periodic with the same period in a neighborhood
of $T$. (See \cite{Radko-classification} for existence of such a
modular vector field). Let now $N$ be an annular neighborhood of
$T$, obtained by translating $N_{0}$ along the flow of $\xi$. We
will extend the identity bisection $\varepsilon_{0}:N_{0}\to X_{0}\cong\Gamma(N_{0})\subset X$
to an identity bisection $\varepsilon:N\to X$. Let $\zeta$ be the
unique lift of $\xi$ to $X$ satisfying\[
(s_{X})_{*}\zeta=(t_{X})_{*}\zeta=\xi.\]
(See \cite{Ginzburg-Lu} for the existence of such a lift). Let $\Phi_{t}$
be the flow of $\xi$, and $\Psi_{t}$ be the flow of $\zeta$ at
time $t$. Note also that on $N_{0}$ the lift of the modular vector
field satisfies\[
\varepsilon_{0}(\Phi_{t}(p))=\Psi_{t}(\varepsilon_{0}(p)).\]
Now define $\varepsilon:N\to X$ by
\begin{enumerate}
\item $\varepsilon|_{N_{0}}=\varepsilon_{0}$.
\item $\varepsilon(p)=\Psi_{t}(\varepsilon(p_{0}))$ for $p=\Phi_{t}(p_{0})\in N$. 
\end{enumerate}
To prove that $\varepsilon$ is defined unambiguously , we need to
check that if $T_{\xi}$ is the period of $\xi$ near the zero curve,
then\[
\Psi_{T_{\xi}}(\varepsilon(p))=\varepsilon(p),\quad\forall p\in N.\]
 Let $\gamma:[0,1]\to P$ be a path such that $\gamma(0)=p$ and $\gamma(1)\in T$.
Then the limits\[
\lim_{t\to1}\varepsilon(\gamma(t)),\qquad\lim_{t\to1}\Psi_{T_{\xi}}(\varepsilon(\gamma(t)))=\Psi_{T_{\xi}}(\lim_{t\to1}\varepsilon(\gamma(t)))\]
both exist. By Lemma \ref{lemma:repelling}, this implies\[
\varepsilon(p)=\Psi_{T_{\xi}}(\varepsilon(p)).\]
Thus $\varepsilon$ is a well-defined identity bisection.

To prove uniqueness, assume that $\varepsilon,\varepsilon':N\to X$
both satisfy\[
s_{X}\circ\varepsilon=t_{X}\circ\varepsilon=\id=t_{X}\circ\varepsilon'=s_{X}\circ\varepsilon'.\]
Then for any $p\in N\setminus T$, consider a path $\gamma:[0,1]\to N$
such that $\gamma(t)\in N\setminus T$ for $t\in[0,1)$ and $\gamma(1)\in T$.
Then the existence of the limits\[
\lim_{t\to1}\varepsilon(\gamma(t)),\qquad\lim_{t\to1}\varepsilon'(\gamma(t))\]
implies that $\varepsilon(p)=\varepsilon'(p)$ by Lemma \ref{lemma:repelling}.
The last corollary implies that $\varepsilon$ is uniquely defined. 
\end{proof}
\begin{thm}
\label{thm:static_Picard}Let $\pi$ be a TSS on a surface $P$ and
$Z\subset P$ be the zero set of the Poisson structure. Then \[
\StatPic(P)\cong\mathcal{M}(P\setminus Z).\]

\end{thm}
\begin{proof}
Let $X\in\StatPic(P)$. Let $L\subset P$ be a symplectic leaf with
the boundary $\partial L=\bigcup_{i=1}^{k}T_{i}$. By Lemma \ref{lem:local-bisections-exist},
for each zero curve $T_{i}\subset\partial L$ there exists a canonical
local identity bisection $\varepsilon_{i}:N_{i}\to X$ in a neighborhood
$N_{i}$ of $T_{i}$. Fix points $p_{1},\ldots,p_{k}$ in the collars
$C_{1}=N_{1}\cap L,\ldots,C_{k}=N_{k}\cap L$ of the boundary curves.
As in section \ref{sec:trivializations}, let $\mathcal{P}(L)$ be
the set of pairs $(X,(x_{1},\dots,x_{k}))$, where $x_{i}\in X_{p_{i}}$
are such that $\Hol_{x_{i}}([T_{i}])=[T_{i}]$. Define a map $\psi_{L}:\ \StatPic(P)\to\mathcal{P}(L)$
by\[
\psi_{L}(X)=(X|_{L},\ (\varepsilon_{1}(p_{1}),\ldots,\varepsilon_{k}(p_{k}))\in\mathcal{P}(L).\]
We claim that the map\[
\psi=\Pi_{L}\psi_{L}:\ \StatPic(P)\to\Pi_{L}\mathcal{P}(L)\cong\Pi_{L}\mathcal{M}(L)=\mathcal{M}(P\setminus Z)\]
is an isomorphism. 

To prove the injectivity, observe that if $X\in\ker\psi_{L}$ for
all $L$, then the local identity bisections can be extended to all
of $P$. Thus \emph{$X$} has a global identity bisection and so is
trivial by Lemma \ref{lem:identityBisectionTrivial}.

We next claim that this map is surjective. Note that any element in
$\Pi_{L}\mathcal{P}(L)\cong\Pi_{L}\mathcal{M}(L)$ can be represented
by a symplectomorphism of $P\setminus Z$, which is identity near
$Z$. Such a symplectomorphism extends to a Poisson isomorphism $\varphi$
of $P$. If we set $X$ to be $\Gamma_{\varphi}(P)$, it is easily
seen that $\psi_{L}(X)$ is exactly the element of $\mathcal{P}(L)$
corresponding via the isomorphism with $\mathcal{M}(L)$ to $\varphi$.
Thus $\psi(X)=\varphi$ and hence our map is onto.
\end{proof}
\begin{cor}
\label{cor:Poiss_to_StatPic}Any bimodule in the static Picard group
can be represented by a Poisson diffeomorphism, i.e., $\forall X\in\StatPic(P)$
$\exists\varphi\in\Poiss(P)$ such that $X\cong\Gamma_{\varphi}(P)$. 
\end{cor}
\begin{proof}
For $X\in\StatPic(P)$ choose a symplectomorphism $\varphi\in\Symp(P\setminus Z)$
representing the corresponding class $\psi(X)\in\mathcal{M}(P\setminus Z)$.
Since $X$ preserves the leaf space pointwise, one can choose $\varphi$
to be trivial near $Z$. Since $\varphi$ preserves the restrictions
of a modular vector field to the zero curves, it can be extended to
a Poisson diffeomorphism $\varphi\in\Poiss(P)$ of the surface. 
\end{proof}
\begin{rem}
\label{rem:stat_Pic_noncompact}Note that (in the case of a compact
surface) $\partial(P\setminus Z)$ , and thus $\mathcal{M}(P\setminus Z)=\pi_{0}(\Diff P\ \textrm{fix\ }Z)=\mathcal{M}(P\ \textrm{fix\ $Z$)}$,
the group of isotopy classes of diffeomorphisms of $P$ fixing a neighborhood
of $Z$ pointwise. If $\pi$ is a TSS on a open surface $P$ such
that $Z\cap\partial P=\emptyset$, the answer for the Picard group
is the same\begin{equation}
\StatPic(P)\cong\mathcal{M}(P\ \textrm{fix\ $Z)$.}\label{eq:stat_Pic_noncompact}\end{equation}

\end{rem}
\begin{example}
\label{exa:cyl_1_static}Let $C\simeq I\times S^{1}$ be a cylinder
with coordinates $(r,\theta)$, where $r\in(-1,1)$ and the Poisson
structure\[
\pi=r\partial_{r}\wedge\partial_{\theta}.\]
The symplectic groupoid of $C$ is given by $\Gamma(C)\simeq C\times\mathbb{R}^{2}$
with coordinates $(r,\theta,p,q)$ and the structure maps\begin{eqnarray*}
s((r,\theta,p,q)) & = & (r,\theta),\\
t((r,\theta,p,q)) & = & (re^{p},\ (\theta+q\cdot r)\mod2\pi),\end{eqnarray*}
and a symplectic structure $\Omega$. The source and target fibers
at all points are isomorphic to $\mathbb{\mathbb{R}}^{2}$, and the
isotropy groups are given by\begin{eqnarray*}
\Gamma_{(0,\theta_{0})}(C) & = & \{(0,\theta_{0},p,q)\ |\ (p,q)\in\mathbb{R}^{2}\}\simeq\mathbb{R}^{2},\\
\Gamma_{(r_{0},\theta_{0})}(C) & = & \{(r_{0},\theta_{0},0,\frac{2\pi k}{r_{0}})\ |\  k\in\mathbb{Z}\}\simeq\mathbb{Z},\quad r_{0}\neq0.\end{eqnarray*}
Away from the zero curve, there are $\mathbb{Z}$ choices of an identity
bisection, given by\[
\sigma_{k}=\{(r,\theta,0,\frac{2\pi k}{r})\ |\  k\in\mathbb{Z}\},\quad r\neq0.\]
However, only one of these {}``almost identity bisections'', $\sigma_{0}$,
extends to the identity bisection $\varepsilon=\sigma_{0}:C\to\Gamma(C)$
on the whole cylinder. This leads to the conclusion that the static
Picard group of the cylinder is trivial, $\StatPic(C)=\{ e\}$. This
corresponds to the general answer (\ref{eq:stat_Pic_noncompact})
as follows: there are two $2$-dimensional symplectic leaves, each
diffeomorphic to a cylinder, with one boundary curve being $T$ and
the other --- a boundary curve of $C$. Since $\mathcal{M}(C\ \textrm{fix $Z$})$
is trivial, it follows that $\StatPic(C)$ is trivial. (By comparison,
the (static) Picard group of a \emph{symplectic} cylinder is $\Pic(S)\cong\Out(\mathbb{Z})=\mathbb{Z}_{2}$). 
\begin{example}
\label{exa:cyl_2_static}Let $C\simeq\mathbb{R}\times S^{1}$ be a
cylinder with coordinates $(r,\theta)$, $r\in(-2,2)$, and the Poisson
structure\[
\pi=(r^{2}-1)\partial_{r}\wedge\partial_{\theta},\]
vanishing linearly on two parallel circles, $r=\pm1$. The symplectic
groupoid $\Gamma(C)$ of this structure is given by $\Gamma(C)=C\times\mathbb{R}^{2}$
with coordinates $(r,\theta,p,q)$. The source and target maps are\begin{eqnarray*}
s((r,\theta,p,q)) & = & (r,\theta),\\
t((r,\theta,p,q)) & = & (\alpha(r,p),\ (\theta+q\cdot(r^{2}-1))\ \mod2\pi),\end{eqnarray*}
where $\alpha(r,p)=\frac{(r+1)+(r-1)e^{2p}}{(r+1)-(r-1)e^{2p}},$
and the symplectic structure is\[
\Omega=-2qrdr\wedge dp+dr\wedge dq-d\theta\wedge dp+(r^{2}-1)dp\wedge dq.\]
Similarly to the previous example, the isotropy at each point away
from the zero curves is isomorphic to $\mathbb{Z},$ while at a point
on a zero curve the isotropy is $\mathbb{R}^{2}$. Let $C_{-1}=\{(r,\theta)\in C\ |r<0\}$
and $C_{1}=\{(r,\theta)\in C\ \ |\  r>0\}$ be disjoint cylindrical
neighborhoods of the zero circles $r=-1$ and $r=1$ respectively,
so that $C=C_{-1}\cup C_{1}\cup\{ r=0\},$ and $\{ r=0\}\in\partial C_{1},\partial C_{-1}$
is the common bounding circle. A bimodule $X\in\Pic(C)$ is obtained
by {}``gluing'' two bimodules $X_{-1}\in\Pic(C_{-1})$ and $X_{1}\in\Pic(C_{1})$.
Since by the previous example the Static Picard groups of $C_{1}$
and $C_{-1}$ are trivial, $X_{-1}\cong\Gamma(C_{-1})$ and $X_{1}\cong\Gamma(C_{1})$.
Each of the groupoids $\Gamma(C_{-1})$ and $\Gamma(C_{1})$ has $\mathbb{Z}$
{}``almost identity bisections'', called $\sigma_{k}^{(-1)}$ and
$\sigma_{k}^{(1)}$ respectively (see previous example). The non-isomorphic
bimodules $X\in\Pic(C)$ arise from various mismatches of these almost
identity bisections (i.e., from gluing the unique identity bisection
$\varepsilon^{(-1)}=\sigma_{0}^{(-1)}$ to all possible $\sigma_{k}^{(1)}$,
$k\in\mathbb{Z}$). Thus, $\StatPic(C)\cong\mathbb{Z}.$ According
to our general formula (\ref{eq:stat_Pic_noncompact}),\[
\StatPic(C)=\mathcal{M}(C\ \textrm{fix$\  Z)$$\cong\mathcal{M}(M_{-1}\ \textrm{fix $\{ r=-1\})\times\mathcal{M}(L)$}\times\mathcal{M}(M_{1}\ \textrm{fix }\ \{ r=1\})$},\]
where $M_{-1}=\{ r<-1\}$, $L=\{-1<r<1\}$ and $M_{1}=\{ r>1\}$ are
the two-dimensional symplectic leaves. Note that $\mathcal{M}(M_{i}\ \textrm{fix$\ \{ r=i\}$})$
is trivial for both $i=\pm1$, and $\mathcal{M}(L)\cong\mathbb{Z}$
is generated by the Dehn twists. 
\end{example}
If we {}``close up'' the cylinder in this example to obtain a Poisson
structure vanishing linearly on two parallel non-separating circles
on the torus, the static Picard group of the resulting structure will
be $\mathbb{Z\times\mathbb{Z}}$. In our picture with almost-identity
bisections each of the copies of $\mathbb{Z}$ corresponds to a possible
{}``mismatch'' of identity bisections in each of the symplectic
cylinders between the zero curves. In our general description in terms
of mapping class groups, this corresponds to Dehn twists in each of
the cylinders between the zero curves. 

More generally, a TSS on a cylinder with $n$ zero curves, which are
separating, has the static Picard group isomorphic to $\mathbb{Z}^{\times(n-1)}$.
A TSS on the torus with $n=2k$ zero curves, which are separating,
is isomorphic to $\mathbb{Z}^{n}$. 
\end{example}

\subsection{The Picard group of a TSS}

To a TSS $(P,\pi)$ we associate a graph $\mathcal{G}(P)$ in the
following way:

\begin{enumerate}
\item a vertex of the graph represents a $2$-dimensional symplectic leaf;
\item an edge represents a common bounding zero curve of two symplectic
leaves;
\item each edge is oriented so that it points toward the vertex for which
the Poisson structure is positive with respect to the orientation
of the surface; 
\end{enumerate}
In addition, we label the graph as follows:

\begin{enumerate}
\item a vertex is labeled by the genus of the corresponding leaf;
\item an edge is labeled by the modular period of $\pi$ around the corresponding
zero curve.
\end{enumerate}
Let $\Aut(\mathcal{G}(P))$$\textrm{ be the group of all automorphism of the graph $\mathcal{G}(P)$ }$
and let $G\subset\Aut(\mathcal{G}(P))$ be its subgroup consisting
of automorphisms preserving the labeling. Since any homeomorphism
of the leaf space gives rise to an automorphism of the graph, we have
a natural group homomorphism $\rho:\Aut(\mathcal{L}(P))\to\Aut(\mathcal{G}(P))$\textcolor{blue}{.}

For a bimodule $X\in\Pic(P)$, let $h_{X}\in\Aut(\mathcal{L}(P))$
be the associated homeomorphism of the leaf space, $h_{X}=t_{X}\circ s_{X}^{-1}$.
Let $j:\Poiss(P)\to\Pic(P)$ be the map $\varphi\mapsto\Gamma_{\varphi}(P)$.
Thus we have the following group homomorphisms:\begin{equation}
\Poiss(P)\stackrel{j}{\to}\Pic(P)\stackrel{h}{\to}\Aut(\mathcal{L}(P))\stackrel{\rho}{\to}\Aut(\mathcal{G}(P)).\label{eq:maps_j_h_rho}\end{equation}

\begin{lem}
The graph automorphism $\rho(h_{X})$ induced by an invertible bimodule
$X\in\Pic(P)$ automatically preserves the labeling. I.e., $\textrm{Im$(\rho\circ h)\subset G$. }$

Moreover, any graph automorphism which preserves the labeling comes
from a Poisson automorphism of the structure, i.e., $\forall\theta\in G\subset\Aut(\mathcal{G}(P))$
$\exists\varphi\in\Poiss(P)$ such that $\rho(h(j(\varphi)))=\theta$.
In other words, the composition $\rho\circ h\circ j$ maps $\Poiss(P,\pi)$
onto $G$.
\end{lem}
\begin{proof}
Let $X\in\Pic(P)$ be a bimodule. Since by \cite{BR} the modular
periods are invariant under Morita equivalence, $\rho(h_{X})\in\Aut(\mathcal{G}(P))$
preserves the labeling of the edges. Since for any $2$-dimensional
leaf $L$ the restriction $X|_{L}$ is a Morita equivalence between
$L$ and $h_{X}(L)$, by a result of Ping Xu (\cite{Xu-ME-Poisson}),
we have $\pi_{1}(L)=\pi_{1}(h_{X}(L))$. Since $h_{X}$ is a homeomorphism,
$L$ and $h_{X}(L)$ have the same number of boundary components.
Thus, $L$ and $h_{X}(L)$ have the same genus. Therefore, the labeling
of the vertices is also preserved.

Let $\theta\in G$ be a graph automorphism preserving the labeling.
Since the graph $\mathcal{G}(P)$ completely encodes the topology
of the decomposition of the surface into the $2$-dimensional symplectic
leaves, there exists a diffeomorphism $\alpha\in\Diff(P)$, sending
leaves to leaves, and inducing $\theta\in G\subset\Aut(\mathcal{G}(P))$.
We may furthermore assume that $\alpha$ preserves the restrictions
of a modular vector field to the zero curves. 

Consider the original Poisson structure $\pi$ and the Poisson structure
$\pi'=\alpha_{*}\pi$ induced by the diffeomorphism $\alpha$. The
zero sets of these two structures are the same, and for any zero curve
$T\in Z(\pi)=Z(\pi')$ the restrictions of the modular vector fields
to $T$ are equal. Thus, there exists a smooth function $f\in C^{\infty}(P)$
such that $\alpha_{*}\pi=f\cdot\pi$, with $f\neq0$. 

We claim that $f>0$. Since $f$ is continuous and nonzero, it is
sufficient to prove that $f>0$ at a point. Let $I$ be a segment
of a common zero curve of $\pi$ and $\pi'$ and $(x,y)$ and $(x',y')$
be the coordinates in a neighborhood $N$ of $I$ such that $\pi=x\partial_{x}\wedge\partial_{y}$
and $\pi'=x'\partial_{x'}\wedge\partial_{y'}$. Note that $x=x'=0$
on $I$. Since the restriction of the modular vector field to $I$
is by assumption preserved by $\alpha$, it follows that $\partial_{y}=\partial_{y'}$
on $I$. 

Let $p\in I$ be a point corresponding to $x=y=0$. It follows that
the Jacobian $J_{\alpha}$ of $\alpha$ at $p$ has the form\[
J_{\alpha}(p)=\left(\begin{array}{cc}
a & 0\\
b & 1\end{array}\right),\qquad a,b\in\mathbb{R}.\]
 We have\[
f(x,y)=\frac{x'}{x}\det J_{\alpha}(x,y).\]
Since $x'=a(y)\cdot x+O(x^{2})$, where $a(y)$ is such that $a(0)=a$,
in a neighborhood of $p$, we get finally that \[
f((0,0))=\lim_{x\to0}\frac{x'}{x}\det J_{\alpha}(x,y)=a^{2}>0.\]
Thus $f>0$. In particular, the symplectic forms corresponding to
$\pi$ and $\pi'$ on the same two-dimensional symplectic leaf have
the same sign (with respect to a chosen symplectic form on the surface).

By applying Moser's argument to each two-dimensional symplectic leaf,
we conclude that $\alpha$ is isotopic to a Poisson diffeomorphism
of $\pi$, i.e., there is a family $\alpha_{t}$, $t\in[0,1]$, of
diffeomorphisms such that 
\begin{itemize}
\item $\alpha_{0}=\alpha$;
\item $\alpha_{1}\in\Poiss(P)$;
\item $\alpha_{t}$ maps leaves to leaves;
\item the restriction of $\alpha_{t}$ to the zero set is equal to a translation
by the flow of a modular vector field; 
\end{itemize}
It follows that each $\alpha_{t}$ induces the same automorphism $\theta\in G$,
and so $\theta=\rho(h(j(\alpha_{1})))$.
\end{proof}
\begin{lem}
For a TSS, the map $j:\Poiss(P)\to\Pic(P)$ is surjective.
\end{lem}
\begin{proof}
Let $X\in\Pic(P)$. By the previous Lemma, the map $\rho\circ h\circ j:\Poiss(P)\to G$
is onto. Thus, there exists a Poisson diffeomorphism $\varphi\in(\rho\circ h\circ j)^{-1}((\rho\circ h)(X))\subset\Poiss(P)$
inducing the same graph automorphism as $X$. By composing $X$ with
$\Gamma_{\varphi^{-1}}(P)$ we may assume that the automorphism $\rho(h_{X})$
of the labeled graph induced by $X$ is the identity map. Next, by
composing $X$ with $\Gamma_{\psi}(P)$ (where the $\psi$ is determined
by the flows of a modular vector field around the zero curves), we
may assume that the leaf space automorphism $h_{X}$ fixes pointwise
the set of zero-dimensional leaves. Thus, $X\in\StatPic(P)$. By Corollary
\ref{cor:Poiss_to_StatPic}, $X$ comes from a Poisson automorphism.
\end{proof}
Thus, for a TSS, the map $j:\Poiss(P)\to\Pic(P)$ is onto. Its kernel
consists of the inner Poisson isomorphisms (i.e., the Poisson isomorphisms
implemented by lagrangian bisections, see \cite{BW1}).

\begin{cor}
For a TSS, $\Pic(P)\cong\Out\Poiss(P)$. 
\end{cor}
\begin{defn}
Let $\mathcal{M}(P,\pi)$ be the group of classes of diffeomorphisms
$\varphi:P\to P$ that map zero curves to zero curves and preserve
the restrictions of the modular vector field to the zero curves, up
to isotopies by diffeomorphisms inducing the identity map on the leaf
space. 
\end{defn}
Note that $\mathcal{M}(P\setminus Z)$ can be considered as a subgroup
of $\mathcal{M}(P,\pi)$ consisting of classes of diffeomorphisms
which preserve the leaf space. 

\begin{thm}
\label{thm:Pic_of_TSS}For a TSS, $\Pic(P)\cong\mathcal{M}(P,\pi)$.
\end{thm}
\begin{rem}
By an argument similar to the one in Lemma \ref{lem:symplectic_repr},
we see that $\mathcal{M}(P,\pi)$ is contractible (by isotopies that
fix neighborhoods of zero curves) to the set of Poisson diffeomorphisms
of $P$.
\end{rem}
\begin{proof}
Using the Remark, represent an element of $\mathcal{M}(P,\pi)$ by
a Poisson automorphism $\varphi$. Define the map $\eta:\mathcal{M}(P,\pi)\to\Pic(P)$
by sending $\varphi$ to the associated bimodule $\Gamma_{\varphi}(P)$.
The restriction of this map to $\mathcal{M}(P\setminus Z)\subset\mathcal{M}(P,\pi)$
is the isomorphism between $\mathcal{M}(P\setminus Z)$ and $\StatPic(P)$.
The map $\eta$ is clearly multiplicative. To check that this map
is well-defined, it is enough to verify this on the class of the identity,
which reduces to the fact that the map from $\mathcal{M}(P\setminus Z)\subset\mathcal{M}(P,\pi)$
to $\StatPic(P)$ is well-defined. By theorem \ref{thm:static_Picard},
this map is surjective. The kernel of $\eta$ clearly lies in $\StatPic(P)$,
and hence is trivial, since the restriction of $\eta$ to $\mathcal{M}(P\setminus Z)\subset\mathcal{M}(P,\pi)$
is an isomorphism.
\end{proof}
\begin{rem}
The answer remains the same for a TSS $\pi$ on an open surface $P$
satisfying $Z\cap\partial P=\emptyset$.
\end{rem}
\begin{example}
\label{exa:cyl_1_Pic}For a TSS on the cylinder with one separating
zero curve (see Example \ref{exa:cyl_1_static}), $\Pic(C)\cong\mathbb{Z}_{2}\times\mathbb{T}^{1}\cong\Out\Poiss(C)$,
where $\mathbb{Z}_{2}$ (isomorphic to the group of automorphisms
of the corresponding graph, which preserve the labeling) is generated
by the (orientation-reversing) Poisson diffeomorphism $\Phi(r,\theta)=(-r,\theta)$,
and the $1$-torus $\mathbb{T}^{1}$ is generated by the flow of a
modular vector field around the zero curve. 
\begin{example}
\label{exa:cyl_2_Pic}For a TSS on the cylinder with two separating
zero curves (see Example \ref{exa:cyl_2_static}), we have $\Pic(C)\cong\mathbb{Z}\times\mathbb{T}\times\mathbb{Z}_{2}$,
where $\mathbb{Z}$ is the Static Picard group (generated by the Dehn
twist of the middle symplectic cylinder), $\mathbb{T}^{2}$ is the
torus generated by rotations of the zero curves, and $\mathbb{Z}_{2}$
(isomorphic to the group of graph automorphisms which preserve the
labeling) corresponds to the flip diffeomorphism $\Phi(r,\theta)=(-r,\theta)$.
Denote by $X_{(k,\phi_{1},\phi_{2},\delta)}\in\Pic(C)$ the bimodule
corresponding to $(k,\phi_{1},\phi_{2},\delta)\in\mathbb{Z}\times\mathbb{T}^{2}\times\mathbb{Z}_{2}$.
Then $X_{(k,0,0,1)}\in\Pic(C)\cong\mathbb{Z}$ are the bimodules in
the Static Picard group. Notice that $X_{(k,0,0,1)}$ can be {}``connected''
to the next static bimodule $X_{((k+1),0,0,1)}$ by a path of bimodules
$X(t)=X_{(k,t,0,1)}$, where $t\in[0,2\pi]$, in the full Picard group.
Indeed, we have $X(0)=X_{(k,0,0,1)}$ and $X(2\pi)=X_{(k,2\pi,0,1)}\cong X_{(k+1,0,0,1)}$. 
\end{example}
\end{example}
\bibliographystyle{amsalpha}

\def\cprime{$'$} \def\cprime{$'$}
\providecommand{\bysame}{\leavevmode\hbox to3em{\hrulefill}\thinspace}

\end{document}